\documentclass[aop,MSNbibl,seceqn,dvips]{arximspdf}

\doi{10.1214/12-AOP829} 
\volume{42}
\issue{3}
\pubyear{2014}
\firstpage{1257}
\lastpage{1284}

\makeatletter
\newcommand{\rrvert}{\vert}
\newcommand{\llvert}{\vert}
\newtheorem{theorem}{Theorem}[section]
\newtheorem{prop}[theorem]{Proposition}
\newtheorem{lemma}[theorem]{Lemma}
\newproclaim{rems}{Remarks}
\newproclaim{rem}{Remark}
\makeatother

\begin{document}
\begin{frontmatter}

\title{Stationary distributions for a class of
generalized~Fleming--Viot processes}
\runtitle{Generalized Fleming--Viot processes}

\begin{aug}
\author{\fnms{Kenji} \snm{Handa}\corref{}\ead[label=e1]{handa@ms.saga-u.ac.jp}}
\runauthor{K. Handa}
\affiliation{Saga University}
\address{
Department of Mathematics\\
Saga University\\
Saga 840-8502 \\
Japan\\
\printead{e1}} 
\end{aug}

\received{\smonth{5} \syear{2012}}
\revised{\smonth{11} \syear{2012}}

%
\begin{abstract}
We identify stationary distributions of generalized Fleming--Viot
processes with jump mechanisms specified by certain beta laws together
with a parameter measure. Each of these distributions is obtained from
normalized stable random measures after a suitable biased
transformation followed by mixing by the law of a Dirichlet random
measure with the same parameter measure. The calculations are based
primarily on the well-known relationship to measure-valued branching
processes with immigration.
\end{abstract}

%

\end{frontmatter}

\section{Introduction}\label{sec1}
In the study of population genetics models, it is of great importance
to identify their stationary distributions. Such identifications
provide us with basic information of possible equilibria of the models
and are needed prior to quantitative discussions on statistical
inference. Since \cite{DK99,Hiraba} and \cite{BL}, theory of
generalized Fleming--Viot processes has served as a new area to be
cultivated and has been developed considerably. (See \cite{BB} for an
exposition.) In view of such progress, it seems that we are in a
position to explore the aforementioned problems for some appropriate
subclass of those models. In this respect, it would be natural to think
of the one-dimensional Wright--Fisher diffusion with mutation as a
prototype. This celebrated process is prescribed by its generator
%
\begin{equation}\label{1.1}
A:=\frac{1}{2}x(1-x)\frac{d^2}{dx^2} +\frac{1}{2}
\bigl[c_1(1-x)-c_2x \bigr]\frac{d}{dx},\qquad x \in[0,1],
\end{equation}
where $c_1$ and $c_2$ are positive constants interpreted as mutation
rates. The stationary distribution is a beta distribution
%
\begin{equation}\label{1.2}
{B}_{c_1,c_2}(dx):= \frac{\Gamma(c_1+c_2)}{\Gamma(c_1)\Gamma(c_2)}
x^{c_1-1}(1-x)^{c_2-1}\,dx,
\end{equation}
where $\Gamma(\cdot)$ is the gamma function. In addition, the process
associated with (\ref{1.1}) admits an infinite-dimensional
generalization known as the Fleming--Viot process with
parent-independent mutation, whose stationary distribution is
identified with the law of a Dirichlet random measure.

In the present paper, we consider a problem of finding a class of
generalized Fleming--Viot processes whose stationary distributions can
be identified. As far as the first term on the right-hand side of
(\ref{1.1}) is concerned, its jump-type version has been discussed in
population genetics as the generator of a model with ``occasional
extreme reproduction''. (See Section~1.2 of \cite{BB} for a
comprehensive account.) We additionally need to look for an appropriate
modification of the second term, which should correspond to a
generalization of the mutation mechanism. With these situations in
mind, our problems can be described as follows.
\begin{longlist}[(II)]
\item[(I)] By modifying both mechanisms of reproduction and mutation,
    find a jump process on $[0,1]$ whose generator extends (\ref{1.1})
    and whose stationary distribution can be identified.

\item[(II)] Establish an analogous generalization for the Fleming--Viot
    process with parent-independent mutation.
\end{longlist}

Since these problems are rather vague, it may be worth showing now the
generator we will believe to give an ``answer'' to~(I). For each
$\alpha\in(0,1)$, define an operator $A_{\alpha}$ by
%
\begin{eqnarray}\label{1.3}
\qquad A_{\alpha}G(x)&=& \int_0^1\frac{B_{1-\alpha,1+\alpha}(du)}{u^2} \bigl[xG
\bigl((1-u)x+u \bigr)\nonumber
\\
&&\hspace*{82pt}{}+(1-x)G \bigl((1-u)x \bigr)-G(x) \bigr]
\nonumber\\[-8pt]\\[-8pt]
&&{} +\int_0^1\frac{B_{1-\alpha,\alpha}(du)}{(\alpha+1)u}\bigl[c_1G \bigl((1-u)x+u \bigr)\nonumber
\\
&&\hspace*{85pt}{} +c_2G \bigl((1-u)x\bigr)-(c_1+c_2)G(x) \bigr],\nonumber
\end{eqnarray}
where $G$ are smooth functions on $[0,1]$. Observe that
$A_{\alpha}G(x)\to AG(x)$ as $\alpha\uparrow1$. It should be noted that
$A_{\alpha}$ is a one-dimensional version of the generator of the
process studied in \cite{BBC} if $c_1=c_2=0$. See also \cite{F} and
\cite{FH}. The reader, however, is cautioned that our notation $\alpha$
is in conflict with that of these papers, in which $\alpha$ plays the
same role as $\alpha+1$ in our notation. (We~adopt such notation in
order for the formulae below to be simpler.) The constant $c_1$ (resp.,
$c_2$) in (\ref{1.3}) can be interpreted as the rate of ``simultaneous
mutation'' from one type to the other type and a proportion $u$ of the
individuals with that type, which are supposed to have the frequency
$1-x$ (resp.,~$x$) in the population, are involved in this ``mutation''
event with\vadjust{\goodbreak} intensity $B_{1-\alpha,\alpha}(du)/((\alpha+1)u)$. [Note
that $(1-u)x+u=x+u(1-x)$.] As will be seen in Proposition~\ref{3.1}
below for more general case, the closure of~(\ref{1.3}) with a suitable
domain generates a Feller semigroup on $C([0,1])$, and our main concern
is the equilibrium state of the associated Markov process. It will be
shown in the forthcoming section that a unique stationary distribution
of the process governed by (\ref{1.3}) is identified with
%
\begin{eqnarray}\label{1.4}
\qquad && P_{\alpha,(c_1,c_2)}(dx)
\nonumber\\[-6pt]\\[-12pt]
&&\qquad:=\Gamma(\alpha+1) \int_0^1{B}_{c_1,c_2}(dy)
E_{\alpha,y} \biggl[(Y_1+Y_2)^{-\alpha};
\frac{Y_1}{Y_1+Y_2}\in dx \biggr],\nonumber
\end{eqnarray}
where $E_{\alpha,y}$ denotes the expectation with respect to
$(Y_1,Y_2)$ with law determined by $\log E_{\alpha,y}[e^{-\lambda_1
Y_1-\lambda_2 Y_2}] =-y\lambda_1^{\alpha}-(1-y)\lambda_2^{\alpha}$
($\lambda_1,\lambda_2\ge0$). Again we see that~(\ref{1.4}) with
$\alpha=1$ reduces to (\ref{1.2}).

One might think that (\ref{1.3}) is one of many possible
generalizations of (\ref{1.1}). In fact it arises naturally in the
following manner. It is well-known \cite{Shiga} that the Fleming--Viot
process with parent-independent mutation can be obtained by way of a
normalization and a random time change from a measure-valued branching
diffusion with immigration. (See also \cite{EM} and~\cite{P}.) An
extension of this significant result was shown in \cite{BBC} for a
class of generalized Fleming--Viot processes, which in the
one-dimensional setting corresponds to (\ref{1.3}) with $c_1=c_2=0$.
Moreover, \cite{BBC} proved that such a jump mechanism is necessary for
a generalized Fleming--Viot process to have the above mentioned link to
a measure-valued branching process with immigration (henceforth
MBI-process). Recently, \cite{FH} showed essentially that the second
term of (\ref{1.3}) is required when we additionally take a
generalization of the mutation mechanism into account. Our argument
will be crucially based on this kind of relationship between the
generalized Fleming--Viot process associated with a natural
generalization of (\ref{1.3}) and a certain ergodic MBI-process. That
relationship can be reformulated as a factorization result on the level
of generators and hence is expected to yield also an explicit
connection between stationary distributions. In principle, the problems
(I)~and~(II) can be considered in a unified way. Nevertheless, we shall
discuss (I)~and~(II) separately. This is mainly because the
factorization identity will turn out to yield a correct answer only for
certain restricted cases and in one dimension one can avoid its use by
taking an analytic approach instead (although this does not reveal
clearly the mathematical structure underlying).\looseness=1

The organization of this paper is as follows. Section~\ref{sec2} is
devoted to derivation of~(\ref{1.4}) by purely analytic argument.
Exploiting the relationship to MBI-processes, we show in
Section~\ref{sec3} that the above mentioned answer to (I) has a natural
generalization which settles (II). The irreversibility of the processes
we consider is discussed in Section~\ref{sec4}.\newpage

\section{The one-dimensional model}\label{sec2}
Let $0<\alpha<1$, $c_1>0$ and $c_2>0$ be given. The purpose of this
section is to show that (\ref{1.4}) is a unique stationary distribution
of the process with generator (\ref{1.3}). Analytically, we shall prove
that a probability measure $P$ on $[0,1]$ satisfying
%
\begin{eqnarray}\label{e2.1}
\int_0^1 A_{\alpha}G(x)P(dx)=0
\nonumber\\[-8pt]\\[-8pt]
\eqntext{\mbox{for all }G(x)=\varphi_n(x):=x^n\mbox{ with }n=1,2,\ldots}
\end{eqnarray}
is uniquely identified with (\ref{1.4}).
Actual starting point of the calculations below is
%
\begin{eqnarray}\label{e2.2}
\int_0^1 A_{\alpha}G(x)P(dx)=0
\nonumber\\[-8pt]\\[-8pt]
\eqntext{\mbox{for all }G(x)=G_t(x):=(1+tx)^{-1}\mbox{ with }t>0.}
\end{eqnarray}
The equivalence of (\ref{e2.1}) and (\ref{e2.2}) is
a consequence of uniform estimates
\[
\bigl|A_{\alpha}\varphi_n(x)\bigr|\le\biggl(1+\frac{c_1+c_2}{\alpha+1}
\biggr)2^n,\qquad n=1,2,\ldots,
\]
which can be shown by observing that
%
\begin{eqnarray}\label{2.3}
&& c_1 \bigl((1-u)x+u \bigr)^n+c_2
\bigl((1-u)x \bigr)^n-(c_1+c_2)x^n
\nonumber
\\
&&\qquad = c_1 \bigl[ \bigl((1-u)x+u \bigr)^n-
\bigl((1-u)x+ux \bigr)^n \bigr]\nonumber
\\
&&\quad\qquad{} +c_2x^n\bigl[(1-u)^n- \bigl((1-u)+u \bigr)^n \bigr]\nonumber
\\
&&\qquad = c_1\sum_{k=1}^n \pmatrix{n
\cr k}(1-u)^{n-k}x^{n-k}u^k
\bigl(1-x^k \bigr)
\nonumber\\\\[-16pt]
&&\quad\qquad{} -c_2x^n\sum _{k=1}^n \pmatrix{n\cr k}(1-u)^{n-k}u^k\nonumber
\\
&&\qquad = \sum_{k=1}^n \pmatrix{n\cr k}(1-u)^{n-k}u^k \bigl[c_1x^{n-k}-(c_1+c_2)x^n\bigr]\nonumber
\\
&&\qquad = u\sum_{k=1}^n \pmatrix{n\cr k}(1-u)^{n-k}u^{k-1} \bigl[c_1x^{n-k}-(c_1+c_2)x^n\bigr]\nonumber
\end{eqnarray}
and in particular
\begin{eqnarray*}
&& x \bigl((1-u)x+u \bigr)^n+(1-x) \bigl((1-u)x\bigr)^n-x^n
\\
&&\qquad =  \sum_{k=2}^n \pmatrix{n\cr
k}(1-u)^{n-k}u^k \bigl(x^{n-k+1}-x^n\bigr)
\\
&&\qquad =  u^2\sum_{k=2}^n \pmatrix{n
\cr k}(1-u)^{n-k}u^{k-2} \bigl(x^{n-k+1}-x^n\bigr).
\end{eqnarray*}
Indeed, these bounds ensure that the function
\[
t \mapsto\int_0^1A_{\alpha}G_t(x)P(dx)
=\sum_{n=1}^{\infty}(-t)^n\int
_0^1A_{\alpha}\varphi_n(x)
P(dx)
\]
is real analytic at least for $-1/2<t<1/2$. We prepare a simple lemma
in order to calculate $A_{\alpha}G_{t}$.

%
\begin{lemma}\label{le2.1}
Assume that $b>0$ and $a+b>0$.
\begin{longlist}[(ii)]
\item[(i)] It holds that for any $\theta_1>0$ and $\theta_2>0$
%
\begin{equation}\label{2.4}
\int_0^1\frac{B_{\theta_1,\theta_2}(du)} {
(au+b)^{\theta_1+\theta_2}}
=(a+b)^{-\theta_1}b^{-\theta_2}.
\end{equation}

\item[(ii)] In addition, suppose that $a'\ne a$ and $a'+b>0$. Then
%
\begin{equation}\label{2.5}
\int_0^1\frac{B_{1-\alpha,1+\alpha}(du)}{(au+b)(a'u+b)} =
\frac{1}{\alpha(a-a')b^{1+\alpha}} \bigl[(a+b)^{\alpha}- \bigl(a'+b
\bigr)^{\alpha} \bigr].
\end{equation}
\end{longlist}
\end{lemma}

Equation (\ref{2.4}) is a one-dimensional version of the formula due to
\cite{CR}, which is sometimes referred to as the Markov--Krein
identity. (See, e.g., \cite{VYT04} or (\ref{3.5}) below.) We will give
a self-contained proof based essentially on the well-known relationship
between beta and gamma laws.

\begin{pf*}{Proof of Lemma~\ref{le2.1}}
The proof of (\ref{2.4}) is simply
done by noting that
\[
(a+b)^{-\theta_1}b^{-\theta_2}= \int_0^{\infty}
\frac{dz_1}{\Gamma(\theta_1)} z_1^{\theta_1-1}e^{-(a+b)z_1} \int
_0^{\infty}\frac{dz_2}{\Gamma(\theta_2)} z_2^{\theta_2-1}e^{-b z_2}
\]
and then by change of variables to $u:=z_1/(z_1+z_2)$, $v:=z_1+z_2$.
The proof of~(\ref{2.5}) can be deduced from (\ref{2.4}) with
$\theta_1=1-\alpha$ and $\theta_2=\alpha$ since
$B_{1-\alpha,1+\alpha}(du) =B_{1-\alpha,\alpha}(du)(1-u)/\alpha$ and
\[
\frac{1-u}{(au+b)(a'u+b)} =\frac{1}{(a-a')b} \biggl(\frac{a+b}{au+b}-
\frac{a'+b}{a'u+b} \biggr).
\]\upqed
\end{pf*}
We proceed to calculate $A_{\alpha}G_t$.

%
\begin{lemma}\label{le2.2}
For any $t>0$ and $x\in[0,1]$,
%
\begin{eqnarray}\label{2.6}
A_{\alpha}G_t(x) &=& t\cdot\frac{(1+t)^{\alpha}-1}{\alpha}\cdot
\frac{x(1-x)}{(1+tx)^{2+\alpha}}
\nonumber\\[-8pt]\\[-8pt]
&&{} -\frac{t}{\alpha+1}\cdot\frac{c_1(1-x)(1+t)^{\alpha-1}-c_2x}{(1+tx)^{1+\alpha}}.\nonumber
\end{eqnarray}
\end{lemma}

\begin{pf}
By straightforward calculations
\begin{eqnarray*}
&& c_1G_t \bigl((1-u)x+u \bigr)+c_2G_t
\bigl((1-u)x \bigr)-(c_1+c_2)G_t(x)
\\
&&\qquad =  -\frac{tu}{1+tx} \biggl[\frac{c_1(1-x)}{1+t(1-u)x+tu} -\frac{c_2x}{1+t(1-u)x} \biggr].
\end{eqnarray*}
Replacing $c_1$ and $c_2$ by $x$ and $1-x$, respectively, we get
\begin{eqnarray*}
&& xG_t \bigl((1-u)x+u \bigr)+(1-x)G_t\bigl((1-u)x \bigr)-G_t(x)
\\
&&\qquad =  \frac{t^2u^2x(1-x)}{1+tx}\cdot\frac{1}{(1+t(1-u)x+tu)(1+t(1-u)x)}.
\end{eqnarray*}
Plugging these equalities into (\ref{1.3}) with $G=G_t$ and then
applying Lemma~\ref{le2.1} yield
%
\begin{eqnarray*}
A_{\alpha}G_t(x) & = & \frac{t^2x(1-x)}{1+tx} \int
_0^1\frac{B_{1-\alpha,1+\alpha}(du)} {
(1+t(1-u)x+tu)(1+t(1-u)x)}
\\
&&{} -\frac{t}{(\alpha+1)(1+tx)} \cdot c_1(1-x)\int_0^1
\frac{B_{1-\alpha,\alpha}(du)}{1+t(1-u)x+tu}
\\
& &{} +\frac{t}{(\alpha+1)(1+tx)} \cdot c_2x\int_0^1
\frac{B_{1-\alpha,\alpha}(du)}{1+t(1-u)x}
\\
& = & \frac{t^2x(1-x)}{1+tx}\cdot\frac{1}{\alpha t (1+tx)^{1+\alpha
}}\cdot\bigl[(1+t)^{\alpha}-1
\bigr]
\\
& &{} -\frac{t}{(\alpha+1)(1+tx)} \biggl[\frac{c_1(1-x)}{(1+t)^{1-\alpha
}(1+tx)^{\alpha}} -\frac{c_2x}{(1+tx)^{\alpha}} \biggr],
\end{eqnarray*}
which equals the right-hand side of (\ref{2.6}).
\end{pf}

Next, we are going to characterize stationary distributions $P$
in terms of
%
\begin{equation}\label{2.7}
S_{\alpha}(t):=\int_0^1
\frac{P(dx)}{(1+tx)^{\alpha}},\qquad t\ge0,
\end{equation}
which is a variant of the generalized Stieltjes transform of
order $\alpha$.

%
\begin{prop}
A probability measure $P$ on $[0,1]$ is a stationary distribution of
the process associated with (\ref{1.3}) if and only if $S_{\alpha}$
defined by (\ref{2.7}) satisfies for all $t>0$
%
\begin{eqnarray}\label{2.8}
&&\frac{(1+t)^{\alpha}-1}{\alpha}(1+t)S_{\alpha}''(t)\nonumber
\\
&&\quad{} + \biggl[ \biggl(c_1+1+\frac{1}{\alpha} \biggr)
\bigl((1+t)^{\alpha}-1 \bigr) +c_1+c_2
\biggr]S_{\alpha}'(t)
\\
&&\quad{} +\alpha c_1(1+t)^{\alpha-1}S_{\alpha}(t)=0.\nonumber
\end{eqnarray}
\end{prop}

\begin{pf}
By virtue of Theorem 9.17 in Chapter~4 of \cite{EK}, $P$ is a
stationary distribution of the process associated with $A_{\alpha}$ if
and only if (\ref{e2.1}) [or (\ref{e2.2})] holds. By Lemma~\ref{le2.2},
(\ref{e2.2}) now reads for all $t>0$
\begin{eqnarray*}
&&  -\frac{(1+t)^{\alpha}-1}{\alpha} \int_0^1
\frac{x(1-x)}{(1+tx)^{2+\alpha}}P(dx)
\\
&&\quad{} +\frac{c_1}{\alpha+1}(1+t)^{\alpha-1} \int_0^1
\frac{1-x}{(1+tx)^{1+\alpha}}P(dx)
\\
&&\quad{}  -\frac{c_2}{\alpha+1}\int_0^1
\frac{x}{(1+tx)^{1+\alpha}}P(dx)=0.
\end{eqnarray*}
This equation becomes (\ref{2.8}) by substituting the equalities
\begin{eqnarray*}
-\int_0^1\frac{x(1-x)}{(1+tx)^{2+\alpha}}P(dx) &=&
\frac{1+t}{\alpha(\alpha+1)}S_{\alpha}''(t) +\frac{1}{\alpha}S_{\alpha}'(t),
\\
\int_0^1\frac{1-x}{(1+tx)^{1+\alpha}}P(dx) &=&
\frac{1+t}{\alpha}S_{\alpha}'(t)+S_{\alpha}(t)
\end{eqnarray*}
and
\[
\int_0^1\frac{x}{(1+tx)^{1+\alpha}}P(dx) =-
\frac{1}{\alpha}S_{\alpha}'(t),
\]
all of which are verified easily.
\end{pf}

We now derive (\ref{1.4}) as the unique stationary
distribution we are looking for.
Recall that for each $y\in(0,1)$ we denote by $E_{\alpha,y}$
the expectation with respect to the two-dimensional random
variable $(Y_1,Y_2)$ with joint law determined by
\[
E_{\alpha,y} \bigl[e^{-\lambda_1 Y_1-\lambda_2 Y_2} \bigr] =e^{-y\lambda
_1^{\alpha}-(1-y)\lambda_2^{\alpha}}, \qquad
\lambda_1,\lambda_2\ge0.
\]
By using
$t^{-\alpha}=\Gamma(\alpha)^{-1}\int_0^{\infty}\,dvv^{\alpha-1}e^{-vt}$
$(t>0)$ and Fubini's theorem, observe that
%
\begin{eqnarray}\label{2.9}
&& E_{\alpha,y} \bigl[(tY_1+Y_2)^{-\alpha} \bigr]\nonumber
\\
&&\qquad = \Gamma(\alpha)^{-1}\int_0^{\infty}\,dvv^{\alpha-1}
\exp\bigl[-y(vt)^{\alpha}-(1-y)v^{\alpha} \bigr]
\\
&&\qquad = \frac{1}{\Gamma(\alpha+1)}\cdot\frac{1}{1+(t^{\alpha}-1)y}\nonumber
\end{eqnarray}
for $t\ge0$. In particular,
$E_{\alpha,y} [(Y_1+Y_2)^{-\alpha} ]=1/\Gamma(\alpha+1)$
and hence
%
\begin{eqnarray}\label{2.10}
&& P_{\alpha,(c_1,c_2)}(dx)
\nonumber\\[-8pt]\\[-8pt]
&&\qquad =\Gamma(\alpha+1)\int_0^1{B}_{c_1,c_2}(dy)
E_{\alpha,y} \biggl[(Y_1+Y_2)^{-\alpha};
\frac{Y_1}{Y_1+Y_2}\in dx \biggr]\nonumber
\end{eqnarray}
defines a probability measure on $[0,1]$.
Although for each $y\in(0,1)$
an expression of the distribution function
\[
[0,1]\ni x\mapsto\Gamma(\alpha+1)E_{\alpha,y} \biggl[(Y_1+Y_2)^{-\alpha};
\frac{Y_1}{Y_1+Y_2}\le x \biggr]
\]
is given as the formula (3.2) in \cite{Y}, that is,
\[
\frac{\sin\alpha\pi}{\pi} \int_0^x
\frac{(1-y)(x-u)^{\alpha-1}u^{\alpha}\,du} {
(1-y)^{2}u^{2\alpha}+y^{2}(1-u)^{2\alpha}
+2y(1-y)u^{\alpha}(1-u)^{\alpha}\cos\alpha\pi},
\]
we do not have any explicit form concerning $P_{\alpha,(c_1,c_2)}$
except the case $c_1+c_2=1$. [See Remark (ii) at the end of this
section.]

The main result of this section is the following.

%
\begin{theorem}\label{th2.4}
The process associated with (\ref{1.3}) has
a unique stationary distribution, which
coincides with $P_{\alpha,(c_1,c_2)}$.
\end{theorem}

\begin{pf}
Notice that the existence of a stationary distribution
follows from compactness of the state space $[0,1]$.
(See, e.g., Remark 9.4 in Chapter~4 of \cite{EK}.)
Let $P$ be an arbitrary stationary distribution of
the process associated with (\ref{1.3})
and $S_{\alpha}$ be defined by (\ref{2.7}).
Put
\[
T_{\alpha}(u)=S_{\alpha} \bigl((1+u)^{1/\alpha}-1 \bigr)
\]
for $u\ge0$.
Setting $t=(1+u)^{1/\alpha}-1$ or $u=(1+t)^{\alpha}-1$,
observe that for $u>0$
\[
T_{\alpha}'(u) =\frac{1}{\alpha}(1+u)^{(1/\alpha)-1}
S_{\alpha}'(t)
\]
and
\begin{eqnarray*}
T_{\alpha}''(u) & = & \frac{1}{\alpha}
\biggl(\frac{1}{\alpha}-1 \biggr) (1+u)^{(1/\alpha)-2}S_{\alpha}'(t)
+ \biggl[\frac{1}{\alpha}(1+u)^{(1/\alpha)-1} \biggr]^2
S_{\alpha}''(t)
\\
& = & \biggl(\frac{1}{\alpha}-1 \biggr) (1+u)^{-1}T_{\alpha}'(u)
+\frac{1}{\alpha^2}(1+u)^{(2/\alpha)-2}S_{\alpha}''(t).
\end{eqnarray*}
Hence, $S_{\alpha}'(t)=\alpha(1+u)^{1-(1/\alpha)}T_{\alpha}'(u)$ and
\[
S_{\alpha}''(t) =\alpha^2(1+u)^{2-(2/\alpha)}
\biggl[T_{\alpha}''(u) - \biggl(
\frac{1}{\alpha}-1 \biggr) (1+u)^{-1}T_{\alpha}'(u)
\biggr].
\]
Also, (\ref{2.8}) can be rewritten as
\begin{eqnarray*}
&& \frac{u}{\alpha}(1+u)^{1/\alpha}S_{\alpha}''(t)
+ \biggl[ \biggl(c_1+1+\frac{1}{\alpha} \biggr)u
+c_1+c_2 \biggr]S_{\alpha}'(t)
\\
&&\quad{}  + \alpha c_1(1+u)^{1-(1/\alpha)}S_{\alpha}(t)=0.
\end{eqnarray*}
From these preliminary observations, it is direct to see that the
equation (\ref{2.8}) is transformed into a hypergeometric equation of
the form
%
\begin{eqnarray}\label{2.11}
\qquad u(1+u)T_{\alpha}''(u) + \bigl[
(c_1+c_2 ) + (c_1+2 )u \bigr]T_{\alpha}'(u)
+c_1T_{\alpha}(u)=0,
\nonumber\\[-10pt]\\[-10pt]
\eqntext{u>0.}
\end{eqnarray}
Clearly, $T_{\alpha}(0)=S_{\alpha}(0)=1$.
In addition,
\[
T_{\alpha}'(0)=S_{\alpha}'(0)/\alpha=-
\int_0^1P(dx)x=-c_1/(c_1+c_2),
\]
where the last equality follows from (\ref{e2.1}) with $n=1$.
These facts together imply that
\[
T_{\alpha}(u) =\int_0^1
\frac{{B}_{c_1,c_2}(dy)}{1+uy}, \qquad u\ge0
\]
or
\[
S_{\alpha}(t) =\int_0^1
\frac{{B}_{c_1,c_2}(dy)} {
1+ \{(1+t)^{\alpha}-1 \}y}, \qquad t\ge0.
\]
(See, e.g., Sections~7.2 and 9.1 in \cite{Lebedev}.)
Combining this with
\begin{eqnarray*}
&& \frac{1}{1+ \{(1+t)^{\alpha}-1 \}y}
\\
&&\qquad =\Gamma(\alpha+1)\int_0^1
\frac{1}{(1+tx)^{\alpha}} E_{\alpha,y} \biggl[(Y_1+Y_2)^{-\alpha};
\frac{Y_1}{Y_1+Y_2}\in dx \biggr],
\end{eqnarray*}
which is immediate from (\ref{2.9}), we arrive at
%
\begin{equation}\label{2.12}
S_{\alpha}(t) =\int_0^1
\frac{P_{\alpha,(c_1,c_2)}(dx)}{(1+tx)^{\alpha}}, \qquad t\ge0
\end{equation}
in view of (\ref{2.10}). Therefore, we conclude that
$P=P_{\alpha,(c_1,c_2)}$ and the proof of Theorem~\ref{th2.4} is
complete.
\end{pf}

\begin{rems*} (i) In the case where $c_1+c_2>1$, an
    alternative  expression for $P_{\alpha,(c_1,c_2)}$ exists:
%
\begin{eqnarray}\label{2.13}
&& P_{\alpha,(c_1,c_2)}(dx)\nonumber
\\
&&\qquad = \Gamma(\alpha+1) (c_1+c_2-1) E
\biggl[(Z_1+Z_2)^{-\alpha}; \frac{Z_1}{Z_1+Z_2}\in dx\biggr]
\\
&&\qquad =: \widetilde{P}_{\alpha,(c_1,c_2)}(dx),\nonumber
\end{eqnarray}
where $Z_1$ and $Z_2$ are independent random variables with
Laplace transforms
%
\begin{equation}\label{2.14}
E \bigl[e^{-\lambda Z_i} \bigr] =\exp\bigl[-c_i\log\bigl(1+
\lambda^{\alpha} \bigr) \bigr], \qquad\lambda\ge0.
\end{equation}
This reflects the fact that the solution to
(\ref{2.11}) with the same initial conditions
$T_{\alpha}(0)=1$ and $T_{\alpha}'(0)=-c_1/(c_1+c_2)$
admits another integral expression of the form
\[
T_{\alpha}(u) =\int_0^1
\frac{{B}_{1,c_1+c_2-1}(dy)} {(1+uy)^{c_1}}, \qquad u\ge0
\]
and accordingly by (\ref{2.12})
%
\begin{equation}
\int_0^1\frac{P_{\alpha,(c_1,c_2)}(dx)}{(1+tx)^{\alpha}} =\int
_0^1\frac{{B}_{1,c_1+c_2-1}(dy)} {
[1+ \{(1+t)^{\alpha}-1 \}y ]^{c_1}}, \qquad t\ge0.
\label{2.15}
\end{equation}
On the other hand, it is not difficult to show that (\ref{2.15}) with
$\widetilde{P}_{\alpha,(c_1,c_2)}$ in place of ${P}_{\alpha,(c_1,c_2)}$
holds, too. In fact, we prove in Lemma~\ref{le3.5} below a
generalization of the coincidence (\ref{2.13}) in the setting of random
measures. Also, the role of $Z_1$~and~$Z_2$ will be made clear in
connection with branching processes with immigration related closely to
the process generated by (\ref{1.3}). [Compare (\ref{2.14}) with
(\ref{3.8}) below.]
\begin{longlist}[(iii)]
\item[(ii)] It will be shown in the Remark after Lemma~\ref{le3.5}
    below that $P_{\alpha,(c_1,c_2)}=B_{\alpha c_1,\alpha c_2}$ holds
    whenever $c_1+c_2=1$. At least at a formal level, this would be
    seen by letting $c_1+c_2\downarrow1$ in (\ref{2.15}) and then by
    making use of (\ref{2.4}).

\item[(iii)] In contrast with the case of the Wright--Fisher diffusion
    mentioned in the \hyperref[sec1]{Introduction}, $P_{\alpha,(c_1,c_2)}$ with
    $0<\alpha<1$ is not a reversible distribution for the generator
    (\ref{1.3}) at least in case $c_1\ne c_2$. This will be seen in
    Section~\ref{sec4}.
\end{longlist}
\end{rems*}

\section{The measure-valued process case}\label{sec3}
The main subject of this section is an extension of Theorem~\ref{th2.4}
to a class of generalized Fleming--Viot processes. But the strategy
will be different from that in the previous section, and so an
alternative proof of Theorem~\ref{th2.4} will be given as a by-product.
To discuss in the setting of measure-valued processes, we need new
notation. Let $E$ be a compact metric space having at least two
distinct points and $C(E)$ [resp., $B_+(E)$] the set of continuous
(resp., nonnegative, bounded Borel) functions on $E$. Define
$\mathcal{M}(E)$ to be the totality of finite Borel measures on $E$,
and we equip $\mathcal{M}(E)$ with the weak topology. Denote by
$\mathcal{M}(E)^{\circ}$ the set of nonnull elements of
$\mathcal{M}(E)$. The set $\mathcal{M}_1(E)$ of Borel probability
measures on $E$ is regarded as a subspace of $\mathcal{M}(E)$. We also
use the notation $\langle\eta, f\rangle$ to stand for the integral of a
function $f$ with respect a measure~$\eta$.
For each $r\in E$, let $\delta_r$ denote
the delta distribution at $r$. Given a probability measure~$Q$,
we write also $E^Q[\cdot]$ for the expectation with
respect to $Q$.

Let $0<\alpha<1$ and $m\in\mathcal{M}(E)$ be given.
We shall discuss in this section
an $\mathcal{M}_1(E)$-valued Markov process associated with
%
\begin{eqnarray}\label{3.1}
&& \mathcal{A}_{\alpha,m}\Phi(\mu)\nonumber
\\[-1pt]
&&\qquad:= \int_0^1 \frac{B_{1-\alpha,1+\alpha}(du)}{u^2}\int
_E\mu(dr) \bigl[\Phi\bigl((1-u)\mu+u\delta_r
\bigr)-\Phi(\mu) \bigr]
\nonumber\\[-8pt]\\[-8pt]
&&\quad\qquad{} +\int_0^1\frac{B_{1-\alpha,\alpha}(du)}{(\alpha+1)u} \int
_E m(dr) \bigl[\Phi\bigl((1-u)\mu+u\delta_r
\bigr)-\Phi(\mu) \bigr],\nonumber
\\[-3pt]
\eqntext{\mu\in\mathcal{M}_1(E),}
\end{eqnarray}
where $\Phi$ belongs to the class $\mathcal{F}_1$ of functions of the
form $\Phi_f(\mu):=\langle\mu^{\otimes n},f \rangle$ for some positive
integer $n$ and $f\in C(E^n)$. Equation~(\ref{3.1}) shows clearly that
$\mathcal{A}_{\alpha,m}$ satisfies the positive maximum principle and
hence is dissipative. (See Lemma~2.1 in Chapter~4 of \cite{EK}.) We
begin by seeing that $\mathcal{A}_{\alpha,m}$ defines a Markov process
on $\mathcal{M}_1(E)$ in an appropriate sense. For this purpose, we
need an expression for $\mathcal{A}_{\alpha,m}\Phi_{f}$ with $f\in
C(E^n)$. Set $(a)_b=\Gamma(a+b)/\Gamma(a)$ for $a>0$ and $b\ge0$, and
let $|\cdot|$ stand for the cardinality. It holds that for any
$\theta\ge0$ and $\nu\in\mathcal{M}_1(E)$
%
\begin{eqnarray}\label{3.1a}
&& \mathcal{A}_{\alpha,\theta\nu}\Phi_f(\mu)\nonumber
\\[-2pt]
&&\qquad  = \sum_{k=2}^n \frac{(1-\alpha)_{k-2} (\alpha+1)_{n-k}}{\Gamma(n)} \sum
_{I\dvtx |I|=k} \bigl( \bigl\langle\mu^{\otimes n},
\Theta^{(n)}_If \bigr\rangle-\Phi_f(\mu)
\bigr)
\\
&&\quad\qquad{} + \theta\sum_{k=1}^n
\frac{(1-\alpha)_{k-1} (\alpha)_{n-k}}{(\alpha
+1)\Gamma(n)} \sum_{I\dvtx |I|=k} \bigl( \bigl\langle
\mu^{\otimes n},\Xi^{(n)}_{I,\nu}f \bigr\rangle-
\Phi_f(\mu) \bigr),\nonumber
\end{eqnarray}
where $I$ are nonempty subsets of $\{1,\ldots,n\}$,
$\Theta^{(n)}_I\dvtx C(E^n)\to C(E^n)$ is defined by letting
$\Theta^{(n)}_If$ be the function obtained from $f$ by replacing all
the variables~$r_i$ with $i\in I$ by $r_{\min I}$ and
$\Xi^{(n)}_{I,\nu}\dvtx C(E^n)\to C(E^n)$ is defined by letting
$\Xi^{(n)}_{I,\nu}f$ be the function obtained from $f$ by replacing all
the variables $r_i$ with $i\in I$ by $r$ and then by integrating with
respect to $\nu(dr)$. (For a degenerate $\nu$, (\ref{3.1a}) is a
special case of the corresponding expression found in the proof of
Lemma 11 in \cite{F}.) Equation~(\ref{3.1a}) can be deduced from the following
identities [cf.~(\ref{2.3})] among signed measures on~$E^n$:
\begin{eqnarray*}
&& \bigotimes_{i=1}^n \bigl((1-u)
\mu(dr_i)+u\delta_r(dr_i) \bigr) -
\bigotimes_{i=1}^n\mu(dr_i)
\\
&&\qquad =  \bigotimes_{i=1}^n \bigl((1-u)
\mu(dr_i)+u\delta_r(dr_i) \bigr) -
\bigotimes_{i=1}^n \bigl((1-u)
\mu(dr_i)+u\mu(dr_i) \bigr)
\\
&&\qquad =  \sum_{I\neq\varnothing} \bigotimes
_{j\notin I} \bigl((1-u)\mu(dr_j) \bigr) \biggl[
\bigotimes_{i\in I} \bigl(u\delta_r(dr_i)
\bigr) -\bigotimes_{i\in I} \bigl(u
\mu(dr_i) \bigr) \biggr]
\\
&&\qquad =  \sum_{I\neq\varnothing} u^{|I|}(1-u)^{n-|I|}
\bigotimes_{j\notin I}\mu(dr_j) \biggl[
\bigotimes_{i\in I}\delta_r(dr_i)
-\bigotimes_{i\in I}\mu(dr_i) \biggr].
\end{eqnarray*}
As for the Fleming--Viot process with parent-independent mutation, the
result corresponding to the next proposition is a special case of
Theorem~3.4 in \cite{EK93}.

%
\begin{prop}\label{pr3.1}
For each $m\in\mathcal{M}(E)$
the closure of $\mathcal{A}_{\alpha,m}$ defined on $\mathcal{F}_1$
generates a Feller semigroup on $C(\mathcal{M}_1(E))$.
\end{prop}

\begin{pf}
Let $\theta\ge0$ and $\nu\in\mathcal{M}_1(E)$ be such that
$m=\theta\nu$. We simply mimic the proof of Theorem~3.4 in \cite{EK93}.
In particular, the Hille--Yosida theorem (Theorem~2.2 in Chapter~4 of
\cite{EK}) will be applied. Let $n$ be an arbitrary positive integer.
Rewrite~(\ref{3.1a}) as
\[
\mathcal{A}_{\alpha,\theta\nu}\Phi_f(\mu) =  \bigl\langle
\mu^{\otimes n},\Theta^{(n)}f \bigr\rangle+\theta\bigl\langle
\mu^{\otimes n},\Xi^{(n)}_{\nu}f \bigr\rangle
-c_n(\alpha,\theta)\Phi_f(\mu),
\]
where $\Theta^{(n)}$, $\Xi^{(n)}_{\nu}\dvtx C(E^n)\to C(E^n)$ and
$c_n(\alpha,\theta)$ are, respectively, the nonnegative operators and
the positive constant-defined implicitly by the above equation combined
with (\ref{3.1a}). Let $\lambda>0$ be arbitrary. Given $g\in C(E^n)$,
define
\[
h= \bigl(\lambda+c_n(\alpha,\theta) \bigr)^{-1}\sum
_{k=0}^{\infty} \bigl[ \bigl(
\lambda+c_n(\alpha,\theta) \bigr)^{-1} \bigl(
\Theta^{(n)}+\theta\Xi^{(n)}_{\nu} \bigr)
\bigr]^k g.
\]
Then $h\in C(E^n)$ since the operator norm of
$\Theta^{(n)}+\theta\Xi^{(n)}_{\nu}$
equals $c_n(\alpha,\theta)$. Moreover,
\[
\bigl(\lambda+c_n(\alpha,\theta) \bigr)h- \bigl(
\Theta^{(n)} +\theta\Xi^{(n)}_{\nu} \bigr)h=g,
\]
so $(\lambda-\mathcal{A}_{\alpha,\theta\nu})\Phi_h=\Phi_g$. This
implies that the range of $\lambda-\mathcal{A}_{\alpha,\theta\nu}$
contains $\mathcal{F}_1$, which is dense in $C(\mathcal{M}_1(E))$. The
rest of the proof is the same as that of Theorem~3.4 in~\cite{EK93}.
\end{pf}

For simplicity, we call the $\mathcal{A}_{\alpha,m}$-process the Markov
process governed by $\mathcal{A}_{\alpha,m}$ in the sense of
Proposition~\ref{pr3.1}. This process is a natural generalization of
the process generated by (\ref{1.3}) in the following sense. Suppose
that $E$ consists of two points, say $r_1$ and $r_2$, set
$m=c_1\delta_{r_1}+c_2\delta_{r_2}$, and let $\{X(t)\dvtx  t\ge0\}$ be
the process generated by (\ref{1.3}). Then, verifying the identity
$\mathcal{A}_{\alpha,m}\Phi(\mu)=A_{\alpha}G(x)$ for
$\mu=x\delta_{r_1}+(1-x)\delta_{r_2}$ and $\Phi(\mu)=G(x)$, we see that
the process $\{X(t)\delta_{r_1}+(1-X(t))\delta_{r_2}\dvtx  t\ge0\}$
defines an $\mathcal{A}_{\alpha,m}$-process. We note that \cite{FH}
discusses the case where $E=[0,1]$ and $m=c\delta_0$ for some $c>0$.

We could also establish the well-posedness of the martingale problem
for $\mathcal{A}_{\alpha,m}$ by modifying some existing arguments. More
precisely, the existence could be shown through a limit theorem for
suitably generalized Moran particle systems by modifying those
considered in the proof of Theorem~2.1 [especially (2.2)]
of~\cite{Hiraba}, which took account of the jump mechanism describing
simultaneous reproduction (sampling) only, so that simultaneous
movement (mutation) of particles to a random location (type)
distributed according to $m(dr)/m(E)$ is allowed. The uniqueness would
follow by the duality argument employing a function-valued process as
in the proof of Theorem~2.1 of \cite{Hiraba}. Its possible transitions
and the associated transition rates are found in (\ref{3.1a}). The
duality would be useful in discussing (weak) ergodicity of the
$\mathcal{A}_{\alpha,m}$-process. (See, e.g., Theorem~5.2 in
\cite{EK93} for such a result in the Fleming--Viot process case.)

The following argument is based primarily on the relationship between
the $\mathcal{A}_{\alpha,m}$-process and a suitable MBI-process, which
takes values in $\mathcal{M}(E)$. More precisely, the generator, say
$\mathcal{L}_{\alpha,m}$, of the latter will be chosen so that for some
constant $C>0$
%
\begin{equation}\label{3.2}
\mathcal{L}_{\alpha,m}\Psi(\eta) =C{\eta(E)^{-\alpha}}
\mathcal{A}_{\alpha,m} \Phi\bigl(\eta(E)^{-1}\eta\bigr), \qquad
\eta\in\mathcal{M}(E)^{\circ},
\end{equation}
where $\Psi(\eta)=\Phi(\eta(E)^{-1}\eta)$ and $\Phi$ is in the linear
span $\mathcal{F}_0$ of functions of the form $\mu\mapsto\langle\mu,
f_1\rangle\cdots\langle\mu, f_n\rangle$ with $f_i\in C(E)$,
$i=1,\ldots,n$ and $n$ being a positive integer. In the case of the
Fleming--Viot process (which corresponds to $\alpha=1$ formally), such
a relation is well known. For instance, it played a key role in
\cite{Shiga}. As for the generalized Fleming--Viot process,
factorizations of the form (\ref{3.2}) have been shown in \cite{BBC}
for $m=0$ (the null measure) and in \cite{FH} for degenerate measures
$m$. From now on, suppose that $m\in\mathcal{M}(E)^{\circ}$. To exploit
(\ref{3.2}) in the study of stationary distributions, we further
require the MBI-process associated with $\mathcal {L}_{\alpha,m}$ to be
ergodic, that is, to have a unique stationary distribution, say
$\widetilde{Q}_{\alpha,m}$, supported on $\mathcal{M}(E)^{\circ}$. Once
these requirements are fulfilled, (\ref{3.2}) suggests that
%
\begin{equation}\label{3.3}
\widetilde{P}_{\alpha,m}(\cdot):= E^{\widetilde{Q}_{\alpha,m}} \bigl[
\eta(E)^{-\alpha}; \eta(E)^{-1}\eta\in\cdot\bigr]
/E^{\widetilde{Q}_{\alpha,m}} \bigl[\eta(E)^{-\alpha} \bigr]
\end{equation}
would give a stationary distribution of the
$\mathcal{A}_{\alpha,m}$-process provided that $\eta(E)^{-\alpha}$ is
integrable with respect to $\widetilde{Q}_{\alpha,m}$. This conditional
answer may be modified to be a general one, which must be consistent
with the one-dimensional result (\ref{1.4}).

To describe the answer, we need both the $\alpha$-stable random measure
with parameter measure $m$ and the Dirichlet random measure with
parameter measure $m$, whose laws on $\mathcal{M}(E)^{\circ}$ and
$\mathcal{M}_1(E)$ are denoted by $Q_{\alpha,m}$ and $\mathcal{D}_{m}$,
respectively. These infinite-dimensional laws are determined uniquely
by the identities
%
\begin{equation}\label{3.4}
\int_{\mathcal{M}(E)^{\circ}}Q_{\alpha,m}(d\eta) e^{-\langle\eta,f
\rangle}
=e^{-\langle m,f^{\alpha} \rangle}
\end{equation}
and
%
\begin{equation}\label{3.5}
\int_{\mathcal{M}_1(E)}\mathcal{D}_{m}(d\mu) {\langle\mu,1+f
\rangle}^{-m(E)} =e^{-\langle m,\log(1+f)\rangle},
\end{equation}
where $f\in B_+(E)$ is arbitrary. A random measure with law
$Q_{\alpha,m}$ is constructed from a Poisson random measure on
$(0,\infty)\times E$. (See also Definition 6 in \cite{VYT04}.) Observe
from (\ref{3.4}) that ${E}^{Q_{\alpha,m}}[\eta(E)^{-\alpha}]
=1/(m(E)\Gamma(\alpha+1))$. As in \cite{Ferguson}, $\mathcal{D}_{m}$~is
defined originally to be the law of a random measure whose arbitrary
finite-dimensional distributions are Dirichlet distributions with
parameters specified by~$m$. The useful identity (\ref{3.5}) is due to
\cite{CR} and reduces to (\ref{2.4}) in one-dimension. We now state the
main result of this paper.

%
\begin{theorem}\label{th3.2}
For any $m\in\mathcal{M}(E)^{\circ}$, the
$\mathcal{A}_{\alpha,m}$-process has a unique stationary distribution,
which is identified with
%
\begin{equation}\label{3.6}
P_{\alpha,m}(\cdot):= \Gamma(\alpha+1)\int_{\mathcal{M}_1(E)}
\mathcal{D}_{m}(d\mu) {E}^{Q_{\alpha,\mu}} \bigl[\eta(E)^{-\alpha};
\eta(E)^{-1}\eta\in\cdot\bigr].
\end{equation}
\end{theorem}

To illustrate, consider the trivial case where $m=\theta\delta_r$ for
some $\theta>0$ and $r\in E$. Then it is verified easily that
$P_{\alpha,m}$ concentrates at $\delta_r\in\mathcal{M}_1(E)$, and this
is consistent with the equality
$\mathcal{A}_{\alpha,m}\Phi(\delta_r)=0$ in that case. Also, for every
$m\in\mathcal{M}(E)^{\circ}$, we note that
$P_{\alpha,m}\to\mathcal{D}_m$ as $\alpha\uparrow1$ since by
(\ref{3.4}) $Q_{\alpha,\mu}$ converges weakly to the delta distribution
at $\mu$ for each $\mu\in\mathcal{M}_1(E)$.

The proof of Theorem~\ref{th3.2} will be divided into three steps. As
mentioned earlier, we first find an ergodic MBI-process whose generator
satisfies (\ref{3.2}) and show, under necessary integrability
condition, that $\widetilde{P}_{\alpha,m}$ in (\ref{3.3}) gives a
stationary distribution of the $\mathcal{A}_{\alpha,m}$-process. [In
fact, the condition will turn out to be that $m(E)>1$. This motivates
us to make a reparametrization $m=:\theta\nu$ with $\theta>0$ and~$\nu\in\mathcal{M}_1(E)$.] Second, for each $\nu\in\mathcal{M}_1(E)$,
we prove that $\widetilde{P}_{\alpha,\theta\nu}=P_{\alpha,\theta\nu}$
for any $\theta>1$. As the last step, we extend stationarity of
$P_{\alpha,\theta\nu}$ with respect to $\mathcal{A}_{\alpha,\theta\nu}$
to all $\theta>0$ by interpreting the condition of stationarity as
certain recursion equations among moment measures which are seen to be
real analytic in $\theta>0$. Also, the recursion equations will be
shown to yield uniqueness of the stationary distribution.

For the first step, we prove in the next proposition
that the MBI-process with the following generator
is the desired one:
%
\begin{eqnarray}\label{3.7}
\qquad&& \mathcal{L}_{\alpha,m}\Psi(\eta)\nonumber
\\
&&\qquad:=  \frac{\alpha+1}{\Gamma(1-\alpha)}\int_0^{\infty}
\frac{dz}{z^{2+\alpha}}\int_E\eta(dr) \biggl[\Psi(\eta+z
\delta_r)-\Psi(\eta)-z\frac{\delta\Psi}{\delta\eta}(r) \biggr]
\nonumber\\[-4pt]\\[-12pt]
&&\quad\qquad{}   -\frac{1}{\alpha} \biggl\langle\eta,\frac{\delta\Psi}{\delta\eta} \biggr\rangle\nonumber
\\
&&\quad\qquad{} +\frac{\alpha}{\Gamma(1-\alpha)}\int_0^{\infty}
\frac{dz}{z^{1+\alpha}}\int_Em(dr) \bigl[\Psi(\eta+z
\delta_r)-\Psi(\eta) \bigr],\nonumber
\end{eqnarray}
where $\Psi$ is in the class $\mathcal{F}$ of functions of the form
$\eta\mapsto F(\langle\eta,f_1\rangle,\ldots,\langle\eta,f_n\rangle )$
for some $F\in C_b^2(\mathbf{R}^n)$, $f_i\in C(E)$ and a positive
integer $n$, and $\frac{\delta\Psi}{\delta\eta}(r)
=\frac{d}{d\varepsilon}\Psi(\eta+\varepsilon\delta_r)\rrvert_{\varepsilon=0}$.
Up to this first order differential term, the operator (\ref{3.7}) for
$E=[0,1]$ and $m=c\delta_0$ with $c>0$ is the same as the one discussed
in Lemma 5.5 of \cite{FH}, in which the factorization (\ref{3.2}) has
been proved. Thus, our main observation in the next proposition is
that, keeping the validity of (\ref{3.2}), such an extra term yields
the ergodicity. Note that the generator (\ref{3.7}) is a special case
of the one discussed in Chapter~9 of \cite{Li}. [See (9.25) combined
with (7.12) there for an expression of the generator.] In particular, a
unique solution to the martingale problem for $\mathcal{L}_{\alpha,m}$
defines an $\mathcal{M}(E)$-valued Markov process, which henceforth we
call the \mbox{$\mathcal{L}_{\alpha,m}$-}process. Intuitively, because of
absence of the ``motion process'', the law of this process is
considered as continuum convolution of the continuous-state branching
process with immigration (CBI-process) studied in \cite{KW}. [See
(\ref{3.10}) below.] In addition, Example 1.1 and Theorem 2.3 in
\cite{KW} concern the one-dimensional version of the
$\mathcal{L}_{\alpha,m}$-process without the drift. The latter proved
that the offspring distribution and the distribution associated with
immigration of the approximating branching processes may have
probability generating functions of the form $s+c(1-s)^{\alpha+1}$ and
$1-d(1-s)^{\alpha}$, respectively.

\begin{prop}\label{pr3.3}
Let $m\in\mathcal{M}(E)^{\circ}$. Then $\mathcal{L}_{\alpha,m}$ in
(\ref{3.7}) and $\mathcal{A}_{\alpha,m}$ in (\ref{3.1}) together
satisfy (\ref{3.2}) with $C=\Gamma(\alpha+2)$ and
$\Psi(\eta)=\Phi(\eta(E)^{-1}\eta)$ for any \mbox{$\Phi\in\mathcal{F}_0$}.
Moreover, the $\mathcal{L}_{\alpha,m}$-process has a unique stationary
distribution $\widetilde{Q}_{\alpha,m}$ with Laplace functional
%
\begin{equation}\label{3.8}
\int_{\mathcal{M}(E)^{\circ}}\widetilde{Q}_{\alpha,m}(d\eta)
e^{-\langle\eta,f \rangle} =e^{-\langle m,\log(1+f^{\alpha}) \rangle},\qquad f\in B_+(E).
\end{equation}
\end{prop}

A random measure with law $\widetilde{Q}_{\alpha,m}$ may be called a
Linnik random measure since it is an infinite-dimensional analogue of
the random variable with law sometimes referred to as a (nonsymmetric)
Linnik distribution, whose Laplace transform appeared already in
(\ref{2.14}). It is obtained by subordinating to an $\alpha$-stable
subordinator by a gamma process. (See, e.g., Example 30.8 in
\cite{Sato}.) Namely, letting $\{Y_{\alpha}(t)\dvtx  t\ge0\}$ and
$\{\gamma(t)\dvtx  t\ge0\}$ be independent L\'evy processes such that
\[
E \bigl[e^{-\lambda Y_{\alpha}(t)} \bigr]=e^{-t\lambda^{\alpha}}\quad\mbox{and}\quad E
\bigl[e^{-\lambda\gamma(t)} \bigr]=e^{-t\log(1+\lambda)},\qquad t,\lambda\ge0,
\]
we have for each $c>0$
\[
E \bigl[e^{-\lambda Y_{\alpha}(\gamma(c))} \bigr] =E \bigl[e^{-\gamma
(c)\lambda^{\alpha}} \bigr]=e^{-c\log(1+\lambda^{\alpha})},
\qquad\lambda\ge0.
\]
The first equality implies that
\[
P \bigl(Y_{\alpha} \bigl(\gamma(c) \bigr)\in\cdot\bigr) =\int
_0^{\infty}P \bigl(\gamma(c)\in dt \bigr)P
\bigl(Y_{\alpha}(t)\in\cdot\bigr).
\]
Equation (\ref{3.8}) clearly shows an analogous structure underlying,
that is,
\[
\widetilde{Q}_{\alpha,m}(\cdot) =\int_{\mathcal{M}(E)^{\circ}}
\mathcal{G}_m(d\eta)Q_{\alpha,\eta
}(\cdot),
\]
where $\mathcal{G}_m$ is the law of the standard gamma process on
$(E,m)$. (See Definition~5 in \cite{VYT04}). It is also obvious from
(\ref{3.8}) that, as $\alpha\uparrow1$, $\widetilde{Q}_{\alpha,m}$
converges to $\mathcal{G}_m$. In addition, one can see that
\[
\lim_{\alpha\uparrow1}\mathcal{L}_{\alpha,m}\Psi(\eta) = \biggl
\langle\eta,\frac{\delta^2\Psi}{\delta\eta^2} \biggr\rangle- \biggl
\langle\eta,
\frac{\delta\Psi}{\delta\eta} \biggr\rangle+ \biggl\langle m,\frac
{\delta\Psi}{\delta\eta} \biggr
\rangle=:\mathcal{L}_{m}\Psi(\eta)
\]
for ``nice'' functions $\Psi$, where
$\frac{\delta^2\Psi}{\delta\eta^2}(r) =\frac{d^2}{d\varepsilon^2}
\Psi(\eta+\varepsilon\delta_r)\rrvert_{\varepsilon=0}$. This is a
special case of the generator of MBI-processes discussed in Section~3
of \cite{Stannat}. It has been proved there that $\mathcal{G}_m$ is a
reversible stationary distribution of the process associated with
$\mathcal{L}_{m}$.

\begin{pf*}{Proof of Proposition~\ref{pr3.3}}
As already remarked, if the term\break
$-\alpha^{-1}\langle\eta,\frac{\delta\Psi}{\delta\eta}\rangle$ in
(\ref{3.7}) would vanish, (\ref{3.2}) can be shown by essentially the
same calculations as in the proof of Lemma 17 in \cite{FH}. [In fact,
the change of variable $z=:\eta(E)u/(1-u)$ in the integrals with
respect to $dz$ in (\ref{3.7}) almost suffices for our purpose.] So,
for the proof of (\ref{3.2}), we only need to observe\vspace*{-2pt} that $\langle
\eta,\frac{\delta\Psi}{\delta\eta}\rangle=0$ for $\Psi$ of the form
$\Psi(\eta)=\Phi(\eta(E)^{-1}\eta)$ with $\Phi\in\mathcal{F}_0$. But
this is readily done by giving a specific form of $\Phi$. Indeed, for
$\Phi(\mu)=\langle\mu, f_1\rangle\cdots\langle\mu, f_n\rangle$ the
function $\Psi$ takes the form $\Psi(\eta)=\langle\eta,
f_1\rangle\cdots\langle\eta, f_n\rangle\langle\eta, 1\rangle^{-n}$,
from which it follows that
\[
\frac{\delta\Psi}{\delta\eta}(r) =\sum_{i=1}^n
\frac{f_i(r)\langle\eta, 1\rangle-\langle\eta,
f_i\rangle} {
\langle\eta, 1\rangle^{n+1}}\prod_{j\ne i}\langle\eta,
f_j\rangle.
\]
After integrating with respect to $\eta(dr)$,
the numerator on the right-hand side vanishes.

The argument regarding ergodicity is based on a well-known formula for
Laplace functionals of transition functions. (See (9.18) in \cite{Li}
for a much more general case than ours.) To write it down, we need only
auxiliary functions called $\Psi$-semigroup \cite{KW} because there is
no ``motion process''. These functions form a one-parameter family
$\{\psi(t,\cdot)\}_{t\ge0}$ of nonnegative functions on $[0,\infty)$
and are determined by the equation
%
\begin{equation}\label{3.9}
\frac{\partial\psi}{\partial t}(t,\lambda) =-\frac{1}{\alpha}\psi
(t,\lambda)^{1+\alpha}
-\frac{1}{\alpha}\psi(t,\lambda), \qquad\psi(0,\lambda)=\lambda
\end{equation}
with $\lambda\ge0$ being arbitrary. An explicit
expression is found in Example 3.1 of \cite{Li}:
\[
\psi(t,\lambda) =\frac{e^{-t/\alpha}\lambda} {
[1+(1-e^{-t})\lambda^{\alpha} ]^{1/\alpha}}.
\]
Let $\{\eta_t\dvtx t\ge0\}$ be an $\mathcal{L}_{\alpha,m}$-process, and
for each $\eta\in\mathcal{M}(E)$ denote by $E_{\eta}$ the expectation
with respect to $\{\eta_t\dvtx t\ge0\}$ starting at $\eta$. Then for
any $f\in B_+(E)$ and $t\ge0$
%
\begin{equation}\label{3.10}
E_{\eta} \bigl[e^{-\langle\eta_t, f\rangle} \bigr] =\exp\biggl
[-\langle\eta,
V_tf\rangle-\int_0^t \bigl
\langle m, (V_sf)^{\alpha} \bigr\rangle \,ds \biggr],
\end{equation}
where $V_tf(r)=\psi(t,f(r))$. As $t\to\infty$
the right-hand side converges to
\[
\exp\biggl[-\int_0^{\infty} \bigl\langle m,
(V_tf)^{\alpha} \bigr\rangle \,dt \biggr] = \exp\bigl[- \bigl
\langle m, \log\bigl(1+f^{\alpha} \bigr) \bigr\rangle\bigr]
\]
since by (\ref{3.9})
\[
\frac{d}{dt}\log\bigl(1+ \bigl(V_tf(r) \bigr)^{\alpha}
\bigr) =- \bigl(V_tf(r) \bigr)^{\alpha}.
\]
This shows the ergodicity required and completes the proof.
\end{pf*}

%
\begin{prop}\label{pr3.4}
Suppose that $m(E)>1$ and let $\widetilde{Q}_{\alpha,m}$ be as in
Proposition~\ref{pr3.3}. Then
\[
E^{\widetilde{Q}_{\alpha,m}} \bigl[\eta(E)^{-\alpha} \bigr] = \bigl
(\Gamma(\alpha+1)
\bigl(m(E)-1 \bigr) \bigr)^{-1}.
\]
Moreover,
%
\begin{equation}
\widetilde{P}_{\alpha,m}(\cdot) = \Gamma(\alpha+1) \bigl(m(E)-1 \bigr)
E^{\widetilde{Q}_{\alpha,m}} \bigl[\eta(E)^{-\alpha}; \eta(E)^{-1}\eta
\in\cdot
\bigr] \label{3.11}
\end{equation}
is a stationary distribution of the $\mathcal{A}_{\alpha,m}$-process.
\end{prop}

\begin{pf}
The first assertion is shown by using
$t^{-\alpha}=\Gamma(\alpha)^{-1}\int_0^{\infty}\,dvv^{\alpha-1}e^{-vt}$
$(t>0)$ and (\ref{3.8}) with $f\equiv v$. Indeed, these equalities
together with Fubini's theorem yield
\begin{eqnarray*}
E^{\widetilde{Q}_{\alpha,m}} \bigl[\eta(E)^{-\alpha} \bigr] &=&  \Gamma(
\alpha)^{-1}\int_0^{\infty}\,dvv^{\alpha-1}
\exp\bigl[-m(E)\log\bigl(1+v^{\alpha} \bigr) \bigr]
\\
&=&  \Gamma(\alpha+1)^{-1}\int_0^{\infty}\,dz
\exp\bigl[-m(E)\log(1+z) \bigr]
\\
&=&  \Gamma(\alpha+1)^{-1} \bigl(m(E)-1 \bigr)^{-1}.
\end{eqnarray*}

As in the one-dimensional case, Theorem 9.17 in Chapter~4 of \cite{EK}
reduces the proof of stationarity of (\ref{3.11}) with respect to
$\mathcal{A}_{\alpha,m}$ to showing that
%
\begin{equation}\label{3.12}
\int_{\mathcal{M}_1(E)} \widetilde{P}_{\alpha,m}(d\mu)
\mathcal{A}_{\alpha,m}\Phi(\mu)=0
\end{equation}
for any $\Phi$ of the form
$\Phi(\mu)=\langle\mu,f_1\rangle\cdots\langle\mu,f_n\rangle$ with
$f_i\in C(E)$ and $n$ being a positive integer. Without any loss of
generality, we can assume that $0\le f_i(x) \le1$ for any $x\in E$ and
$i=1,\ldots,n$. Furthermore, we only have to consider the case where
$f_1=\cdots=f_n=:f$ because the coefficients of the monomial $t_1\cdots
t_n$ in $\langle\mu,t_1f_1+\cdots+t_nf_n\rangle^n$ equals
$n!\langle\mu,f_1\rangle\cdots\langle\mu,f_n\rangle$. Thus, we let
$\Phi(\mu)=\langle\mu,f\rangle^n$ with $0\le f(x) \le1$ for any $x\in
E$. Because of the basic relation (\ref{3.2}) and (\ref{3.11})
together, (\ref{3.12}) can be rewritten as
%
\begin{equation}\label{3.13}
\int_{\mathcal{M}(E)^{\circ}} \widetilde{Q}_{\alpha,m}(d\eta)
\mathcal{L}_{\alpha,m}\Psi(\eta)=0,
\end{equation}
where $\Psi(\eta)=\langle\eta,f\rangle^n\langle\eta,1\rangle^{-n}$. The
main difficulty comes from the fact that $\Psi$ does not belong to
$\mathcal{F}$. For each $\varepsilon>0$, introduce
$\Psi_{\varepsilon}(\eta):= \langle\eta,f\rangle^n
(\langle\eta,1\rangle+\varepsilon)^{-n}$ and observe that
$\Psi_{\varepsilon}\in\mathcal{F}$. Thanks to Proposition~\ref{pr3.3},
we then have (\ref{3.13}) with $\Psi_{\varepsilon}$ in place of $\Psi$
provided that $\mathcal{L}_{\alpha,m}\Psi_{\varepsilon}$ is bounded.
Thus, the proof of (\ref{3.13}) reduces to showing the following two
assertions:
\begin{longlist}[(ii)]
\item[(i)] For every $\varepsilon>0$,
    $\mathcal{L}_{\alpha,m}^{(1)}\Psi_{\varepsilon}$,
    $\mathcal{L}_{\alpha,m}^{(2)}\Psi_{\varepsilon}$ and
    $\mathcal{L}_{\alpha,m}^{(3)}\Psi_{\varepsilon}$ are bounded
    functions on~$\mathcal{M}(E)$.

\item[(ii)] It holds that for each $k\in\{1,2,3\}$
%
\begin{equation}
\quad \lim_{\varepsilon\downarrow0}\int_{\mathcal{M}(E)^{\circ}}
\widetilde{Q}_{\alpha,m}(d\eta)\mathcal{L}_{\alpha,m}^{(k)}\Psi
_{\varepsilon}(\eta) = \int_{\mathcal{M}(E)^{\circ}} \widetilde{Q}_{\alpha,m}(d
\eta)\mathcal{L}_{\alpha,m}^{(k)}\Psi(\eta). \label{3.14}
\end{equation}
\end{longlist}
Here, $\mathcal{L}_{\alpha,m}
=\mathcal{L}_{\alpha,m}^{(1)}+\mathcal{L}_{\alpha,m}^{(2)}+\mathcal
{L}_{\alpha,m}^{(3)}$,
and the operators $\mathcal{L}_{\alpha,m}^{(1)}$,
$\mathcal{L}_{\alpha,m}^{(2)}$ and $\mathcal{L}_{\alpha,m}^{(3)}$
correspond, respectively, to the first, second and last term on the
right-hand side of~(\ref{3.7}).

First, we consider $\mathcal{L}_{\alpha,m}^{(2)}$.
Observe that
%
\begin{eqnarray}\label{3.15}
\frac{\delta\Psi_{\varepsilon}}{\delta\eta}(r) &=&  \frac{nf(r)\langle\eta,f\rangle^{n-1}}{
(\langle\eta,1\rangle+\varepsilon)^n} -\frac{n\langle\eta,f\rangle^n} {(\langle\eta,1\rangle+\varepsilon)^{n+1}}
\nonumber\\[-8pt]\\[-8pt]
&=&  \frac{n(f(r)\langle\eta,1\rangle-\langle\eta,f\rangle+\varepsilon
f(r) )\langle\eta,f\rangle^{n-1}} {
(\langle\eta,1\rangle+\varepsilon)^{n+1}},\nonumber
\end{eqnarray}
from which it follows that
\begin{eqnarray*}
\alpha\mathcal{L}_{\alpha,m}^{(2)}\Psi_{\varepsilon}(\eta) &=&
- \biggl\langle\eta,\frac{\delta\Psi_{\varepsilon}}{\delta\eta} \biggr
\rangle
\\
&=&  -\frac{n (\langle\eta,f\rangle\langle\eta,1\rangle -\langle\eta,f\rangle\langle\eta,1\rangle+\varepsilon\langle
\eta,f\rangle ) \langle\eta,f\rangle^{n-1}} {(\langle\eta,1\rangle+\varepsilon)^{n+1}}
\\
&=&  -n\varepsilon\frac{\Psi_{\varepsilon}(\eta)}{\langle\eta,1\rangle
+\varepsilon}.
\end{eqnarray*}
Hence, $\mathcal{L}_{\alpha,m}^{(2)}\Psi_{\varepsilon}$ is a bounded
function on $\mathcal{M}(E)$ and
$\mathcal{L}_{\alpha,m}^{(2)}\Psi_{\varepsilon}(\eta)\to
0=\mathcal{L}_{\alpha,m}^{(2)}\Psi(\eta)$ boundedly as
$\varepsilon\downarrow0$. This proves that (i)~and~(ii) hold true for~$\mathcal{L}_{\alpha,m}^{(2)}$.

In calculating $\mathcal{L}_{\alpha,m}^{(3)}\Psi_{\varepsilon}$,
(\ref{3.15}) is useful since
$\frac{d}{dz}\Psi_{\varepsilon}(\eta+z\delta_r)=
\frac{\delta\Psi_{\varepsilon}}{\delta(\eta+z\delta_r)}(r)$. Indeed, by
Fubini's theorem
%
\begin{eqnarray}\label{3.16}
&& \int_0^{\infty} \frac{dz}{z^{1+\alpha}} \bigl[
\Psi_{\varepsilon}(\eta+z\delta_r)-\Psi_{\varepsilon}(\eta)
\bigr]\nonumber
\\
&&\qquad =  \int_0^{\infty} \frac{dz}{z^{1+\alpha}} \int
_0^z\,dw\frac{\delta\Psi_{\varepsilon}}{\delta(\eta+w\delta_r)}(r)
\\
&&\qquad =  \frac{1}{\alpha}\int_0^{\infty}
w^{-\alpha}\,dw\frac{\delta\Psi_{\varepsilon}}{\delta(\eta+w\delta_r)}(r)\nonumber
\end{eqnarray}
and combining with (\ref{3.15}) yields
%
\begin{eqnarray}\label{3.17}
&&{\fontsize{10pt}{12pt}\selectfont{\mbox{$\displaystyle{\biggl\llvert\int_0^{\infty}
\frac{dz}{z^{1+\alpha}} \bigl[\Psi_{\varepsilon}(\eta+z\delta_r)-
\Psi_{\varepsilon}(\eta) \bigr] \biggr\rrvert}$}}}\nonumber
\\
&&\hspace*{-3pt}{\fontsize{10pt}{12pt}\selectfont{\mbox{$\displaystyle{\qquad \le \frac{1}{\alpha}\int_0^{\infty}w^{-\alpha}\,dw
\frac{n \llvert f(r)\langle\eta+w\delta_r,1\rangle
-\langle\eta+w\delta_r,f\rangle+\varepsilon f(r)\rrvert
\langle\eta+w\delta_r,f\rangle^{n-1}} {(\langle\eta+w\delta_r,1\rangle+\varepsilon)^{n+1}}}$}}}
\nonumber\hspace*{-21pt}
\\
&&\hspace*{-3pt}{\fontsize{10pt}{12pt}\selectfont{\mbox{$\displaystyle{\qquad \le \frac{n}{\alpha}\int_0^{\infty}w^{-\alpha}\,dw
\frac{1}{\langle\eta,1\rangle+w+\varepsilon}}$}}}
\\
&&\hspace*{-3pt}{\fontsize{10pt}{12pt}\selectfont{\mbox{$\displaystyle{\qquad =  \frac{n}{\alpha}\int_0^{\infty}w^{-\alpha}\,dw
\int_0^{\infty}\,dv e^{-v(\langle\eta,1\rangle+w+\varepsilon)}}$}}}
\nonumber
\\
&&\hspace*{-3pt}{\fontsize{10pt}{12pt}\selectfont{\mbox{$\displaystyle{\qquad =  n\frac{\Gamma(\alpha)\Gamma(1-\alpha)}{\alpha} \bigl(\langle\eta,1\rangle+\varepsilon
\bigr)^{-\alpha}.}$}}}\nonumber
\end{eqnarray}
This shows not only that
$\mathcal{L}_{\alpha,m}^{(3)}\Psi_{\varepsilon}$ is bounded but also
\[
\bigl|\mathcal{L}_{\alpha,m}^{(3)}\Psi_{\varepsilon}(\eta)\bigr| \le n
\Gamma(\alpha) \cdot\frac{\langle m,1\rangle}{\langle\eta,1\rangle
^{\alpha}},
\]
which is integrable with respect to $\widetilde{Q}_{\alpha,m}$ as
proved already. It can be seen also from (\ref{3.15}) and (\ref{3.16})
that $\mathcal{L}_{\alpha,m}^{(3)}\Psi_{\varepsilon}$ converges
pointwise to $\mathcal{L}_{\alpha,m}^{(3)}\Psi$ as
$\varepsilon\downarrow0$. By Lebesgue's dominated convergence theorem
we have proved (\ref{3.14}) for $\mathcal{L}_{\alpha,m}^{(3)}$.

The final task is to deal with $\mathcal{L}_{\alpha,m}^{(1)}\Psi
_{\varepsilon}$. Similar to (\ref{3.16})
\begin{eqnarray*}
I_{\varepsilon}(\eta,r) &:= & \int_0^{\infty}
\frac{dz}{z^{2+\alpha}} \biggl[\Psi_{\varepsilon}(\eta+z\delta_r)-
\Psi_{\varepsilon}(\eta) -z\frac{\delta\Psi_{\varepsilon}}{\delta\eta}(r)
\biggr]
\\
&=&  \int_0^{\infty} \frac{dz}{z^{2+\alpha}} \int
_0^z\,dw \biggl[\frac{\delta\Psi_{\varepsilon}}{\delta(\eta+w\delta_r)}(r) -
\frac{\delta\Psi_{\varepsilon}}{\delta\eta}(r) \biggr]
\\
&=&  \frac{1}{1+\alpha}\int_0^{\infty}
\frac{dw}{w^{1+\alpha}} \biggl[\frac{\delta\Psi_{\varepsilon}}{\delta(\eta
+w\delta_r)}(r) -\frac{\delta\Psi_{\varepsilon}}{\delta\eta}(r) \biggr].
\end{eqnarray*}
By (\ref{3.15})
$\frac{\delta\Psi_{\varepsilon}}{\delta(\eta+w\delta_r)}(r)
-\frac{\delta\Psi_{\varepsilon}}{\delta\eta}(r)$ equals
{\fontsize{10.3pt}{12pt}\selectfont{\begin{eqnarray*}
\hspace*{-4pt}&& \frac{(\langle\eta,1\rangle+\varepsilon)^{n+1}n
(f(r)\langle\eta,1\rangle-\langle\eta,f\rangle+\varepsilon
f(r) )[\langle\eta+w\delta_r,f\rangle^{n-1}
-\langle\eta,f\rangle^{n-1} ]} {
(\langle\eta,1\rangle+w+\varepsilon)^{n+1}
(\langle\eta,1\rangle+\varepsilon)^{n+1}}
\\
\hspace*{-4pt}&&\quad{} +  \frac{[(\langle\eta,1\rangle+\varepsilon)^{n+1}
-(\langle\eta,1\rangle+w+\varepsilon)^{n+1} ]
n (f(r)\langle\eta,1\rangle-\langle\eta,f\rangle
+\varepsilon f(r) )\langle\eta,f\rangle^{n-1}} {
(\langle\eta,1\rangle+w+\varepsilon)^{n+1}
(\langle\eta,1\rangle+\varepsilon)^{n+1}}.
\end{eqnarray*}}}%
Moreover, we have bounds
\begin{eqnarray*}
\bigl\llvert\langle\eta+w\delta_r,f\rangle^{n-1} -
\langle\eta,f\rangle^{n-1} \bigr\rrvert &=&  \biggl\llvert\int
_0^w\,dv(n-1)f(r) \langle\eta+v \delta_r,f\rangle^{n-2} \biggr\rrvert
\\
&\le& w(n-1) \bigl(\langle\eta,1\rangle+w \bigr)^{n-2}
\end{eqnarray*}
and
\begin{eqnarray*}
\bigl\llvert\bigl(\langle\eta,1\rangle+\varepsilon\bigr)^{n+1} -
\bigl( \langle\eta,1\rangle+w+\varepsilon\bigr)^{n+1} \bigr\rrvert
&=& (n+1)\int_0^w\,dv \bigl(\langle\eta,1\rangle
+v+\varepsilon\bigr)^{n}
\\
&\le& w(n+1) \bigl(\langle\eta,1\rangle+w+\varepsilon\bigr)^{n}.
\end{eqnarray*}
Consequently,
\begin{eqnarray*}
&& \biggl\llvert\frac{\delta\Psi_{\varepsilon}}{\delta(\eta+w\delta_r)}(r)
-\frac{\delta\Psi_{\varepsilon}}{\delta\eta}(r) \biggr\rrvert
\\
&&\qquad \le w \frac{n(\langle\eta,1\rangle+\varepsilon)^{n+2}
(n-1)(\langle\eta,1\rangle+w)^{n-2}} {
(\langle\eta,1\rangle+w+\varepsilon)^{n+1}
(\langle\eta,1\rangle+\varepsilon)^{n+1}}
\\
&&\quad\qquad{} +w\frac{(n+1)(\langle\eta,1\rangle+w+\varepsilon)^{n}
n(\langle\eta,1\rangle+\varepsilon)\langle\eta,1\rangle^{n-1}} {
(\langle\eta,1\rangle+w+\varepsilon)^{n+1}
(\langle\eta,1\rangle+\varepsilon)^{n+1}}
\\
&&\qquad \le w\frac{2n^2}{(\langle\eta,1\rangle+w+\varepsilon)
(\langle\eta,1\rangle+\varepsilon)}.
\end{eqnarray*}
Therefore, analogous calculations to those in (\ref{3.17}) lead to
\begin{eqnarray*}
\bigl\llvert\mathcal{L}_{\alpha,m}^{(1)}\Psi_{\varepsilon}(
\eta) \bigr\rrvert &=&  \biggl\llvert\frac{\alpha+1}{\Gamma(1-\alpha
)} \int _E I_{\varepsilon}(\eta,r)\eta(dr) \biggr\rrvert
\\
&\le& 2n^2\Gamma(\alpha) \bigl(\langle\eta,1\rangle+\varepsilon
\bigr)^{-\alpha} \cdot\frac{\langle\eta,1\rangle}{\langle\eta,1\rangle +\varepsilon}.
\end{eqnarray*}
This makes it possible to argue as in the case of
$\mathcal{L}_{\alpha,m}^{(3)}\Psi_{\varepsilon}$ to verify (i) and (ii)
for~$\mathcal{L}_{\alpha,m}^{(1)}$. We complete the proof of
Proposition~\ref{pr3.4}.
\end{pf}

Next, we show the coincidence of two distributions
(\ref{3.3}) [or (\ref{3.11})] and (\ref{3.6}).
Before going to the proof, it is worth noting that
%
\begin{equation}\label{3.18}
{P}_{\alpha,m}(\cdot) =\int_{\mathcal{M}_1(E)}\mathcal{D}_{m}(d
\mu)\mathcal{D}^{(\alpha,\alpha)}_{\mu}(\cdot),
\end{equation}
where in general, for $\theta>-\alpha$ and $m\in\mathcal{M}(E)$,
$\mathcal{D}^{(\alpha,\theta)}_{m}$ is the law of the two-parameter
generalization of the Dirichlet random measure with parameter
$(\alpha,\theta)$ and parameter measure $m$ defined by
\[
\mathcal{D}^{(\alpha,\theta)}_{m}(\cdot) = \frac{\Gamma(\theta
+1)}{\Gamma((\theta/\alpha)+1)}
{E}^{Q_{\alpha,m}} \bigl[\eta(E)^{-\theta}; \eta(E)^{-1}\eta\in
\cdot\bigr].
\]
(See, e.g., Section~5 of \cite{VYT04}.)
We will make use of the identity
%
\begin{equation}\label{3.19}
\qquad\qquad \int_{\mathcal{M}_1(E)}\mathcal{D}^{(\alpha,\alpha)}_{m}(d\mu)
\langle\mu, 1+f\rangle^{-\alpha} = \bigl\langle m, (1+f)^{\alpha}
\bigr\rangle^{-1},\qquad f\in B_+(E).
\end{equation}
This is a special case of Theorem 4 in \cite{VYT04} and
can be shown as follows:
\begin{eqnarray*}
&& \int_{\mathcal{M}_1(E)}\mathcal{D}^{(\alpha,\alpha)}_{m}(d\mu)
\langle\mu, 1+f\rangle^{-\alpha}
\\
&&\qquad =  \Gamma(\alpha+1)E^{{Q}_{\alpha,m}}
\bigl[\langle\eta,1\rangle^{-\alpha} \bigl(1+\langle\eta,1\rangle
^{-1}\langle\eta,f\rangle\bigr)^{-\alpha} \bigr]
\\
&&\qquad =  \Gamma(\alpha+1)E^{{Q}_{\alpha,m}} \bigl[\langle\eta,1+f\rangle
^{-\alpha} \bigr]
\\
&&\qquad =  \alpha\int_0^{\infty}\,dvv^{\alpha-1} \exp
\bigl[-v^{\alpha} \bigl\langle m,(1+f)^{\alpha} \bigr\rangle\bigr]
\\
&&\qquad =  \bigl\langle m, (1+f)^{\alpha} \bigr\rangle^{-1}.
\end{eqnarray*}

%
\begin{lemma}\label{le3.5}
If $m(E)>1$, then $\widetilde{P}_{\alpha,m}$ in (\ref{3.11})
coincides with ${P}_{\alpha,m}$ in~(\ref{3.6}).
\end{lemma}

\begin{pf}
It suffices to show that for any $f\in B_+(E)$
\begin{eqnarray*}
\widetilde{I}(f)&:=&\int_{\mathcal{M}_1(E)}\widetilde{P}_{\alpha,m}(d
\mu) \langle\mu,1+f\rangle^{-\alpha}
\\
&=& \int_{\mathcal{M}_1(E)}{P}_{\alpha,m}(d\mu) \langle\mu,1+f\rangle^{-\alpha}=:I(f).
\end{eqnarray*}
In view of (\ref{3.11}), calculations similar to the proof of
(\ref{3.19}) show that
\begin{eqnarray*}
&& \bigl(\Gamma(\alpha+1) \bigl(m(E)-1 \bigr)
\bigr)^{-1}\widetilde{I}(f)
\\
&&\qquad = E^{\widetilde{Q}_{\alpha,m}} \bigl [\langle\eta,1+f \rangle^{-\alpha}
\bigr]
\\
&&\qquad =  \Gamma(\alpha)^{-1}\int_0^{\infty}\,dvv^{\alpha-1}
\exp\bigl[- \bigl\langle m,\log\bigl(1+v^{\alpha}(1+f)^{\alpha}
\bigr) \bigr\rangle\bigr]
\\
&&\qquad =  \Gamma(\alpha+1)^{-1}\int_0^{\infty}\,dz
\exp\bigl[- \bigl\langle m,\log\bigl(1+z(1+f)^{\alpha} \bigr) \bigr
\rangle\bigr]
\\
&&\qquad =  \frac{1}{\Gamma(\alpha+1)}\int_0^{1}\,du(1-u)^{-2}
\exp\biggl[- \biggl\langle m,\log\biggl(1+\frac{u}{1-u}(1+f)^{\alpha}
\biggr) \biggr\rangle\biggr]
\\
&&\qquad =  \frac{1}{\Gamma(\alpha+1)}\int_0^{1}\,du(1-u)^{m(E)-2}
\exp\bigl[- \bigl\langle m,\log\bigl(1+u \bigl((1+f)^{\alpha}-1 \bigr)
\bigr) \bigr\rangle\bigr]
\\
&&\qquad =  \frac{1}{\Gamma(\alpha+1)}\int_0^{1}\,du(1-u)^{m(E)-2}
\\
&&\quad\qquad{}\times \int_{\mathcal{M}_1(E)}\mathcal{D}_{m}(d\mu) \bigl\langle\mu,
1+u \bigl((1+f)^{\alpha}-1 \bigr)\bigr\rangle^{-m(E)},
\end{eqnarray*}
where the last equality follows from (\ref{3.5}).
Hence, by applying Fubini's theorem and (\ref{2.4})
\begin{eqnarray*}
\widetilde{I}(f) &=&  \int_{\mathcal{M}_1(E)}\mathcal{D}_{m}(d\mu) \int
_0^1\frac{B_{1,m(E)-1}(du)} {
\langle\mu, 1+u((1+f)^{\alpha}-1)\rangle^{m(E)}}
\\
&=&  \int_{\mathcal{M}_1(E)}\mathcal{D}_{m}(d\mu) \bigl
\langle\mu, (1+f)^{\alpha} \bigr\rangle^{-1}.
\end{eqnarray*}
On the other hand, combining (\ref{3.18}) with (\ref{3.19}),
we get
%
\begin{equation}\label{3.20}
I(f)=\int_{\mathcal{M}_1(E)}\mathcal{D}_{m}(d\mu) \bigl
\langle\mu, (1+f)^{\alpha} \bigr\rangle^{-1}
\end{equation}
and therefore $I(f)=\widetilde{I}(f)$ as desired.
\end{pf}

\begin{rem*}
The ``semi-explicit'' form (\ref{3.18}) can be explicit if $m$ is a
probability measure. More precisely, we have
$P_{\alpha,\nu}=\mathcal{D}_{\alpha\nu}$ for any $\nu\in
\mathcal{M}_1(E)$. Indeed, observe that by (\ref{3.20}) with $m=\nu$
\begin{eqnarray*}
\int_{\mathcal{M}_1(E)}{P}_{\alpha,\nu}(d\mu)\langle\mu, 1+f\rangle
^{-\alpha} &=&  \int_{\mathcal{M}_1(E)}\mathcal{D}_{\nu}(d
\mu) \bigl\langle\mu, (1+f)^{\alpha} \bigr\rangle^{-1}
\\
&=&  \exp\bigl[- \bigl\langle\nu, \log\bigl\{(1+f)^{\alpha} \bigr\}
\bigr\rangle\bigr]
\\
&=&  \exp\bigl[- \bigl\langle\alpha\nu, \log(1+f) \bigr\rangle\bigr]
\\
&=&  \int_{\mathcal{M}_1(E)}{\mathcal{D}}_{\alpha\nu}(d\mu)\langle
\mu, 1+f\rangle^{-\alpha},
\end{eqnarray*}
where (\ref{3.5}) has been applied twice. [A one-dimensional version of
the identity $P_{\alpha,\nu}=\mathcal {D}_{\alpha\nu}$ is mentioned in
Remark (ii) at the end of Section~\ref{sec2}.] By (\ref{3.18}) what we
have just seen is rewritten as
\[
\int_{\mathcal{M}_1(E)}\mathcal{D}_{\nu}(d\mu)\mathcal
{D}^{(\alpha,\alpha)}_{\mu}(\cdot) = \mathcal{D}_{\alpha\nu}(\cdot),
\]
which is a special case of
\[
\int_{\mathcal{M}_1(E)}\mathcal{D}^{(\beta,\theta/\alpha)}_{\nu
}(d\mu)
\mathcal{D}^{(\alpha,\theta)}_{\mu}(\cdot) = \mathcal{D}^{(\alpha\beta,\theta)}_{\nu}(\cdot), \qquad\beta\in[0,1), \theta>-\alpha\beta.
\]
Here notice that, in case $\beta=0$,
$\mathcal{D}^{(0,\theta)}_{\nu}=\mathcal{D}_{\theta\nu}$ by definition.
This generalization can be proved analogously by virtue of the
two-parameter generalization of (\ref{3.5}) and (\ref{3.19}). (See,
e.g., Theorem 4 in \cite{VYT04}.)
\end{rem*}

We can now prove our main result, Theorem~\ref{th3.2}. In the proof, we
write $\theta\nu$ [$\theta>0$, $\nu\in\mathcal{M}_1(E)$] for the
parameter measure $m$.

\begin{pf*}{Proof of Theorem~\ref{th3.2}}
Let $\nu\in\mathcal{M}_1(E)$ be
given. We first show that, for arbitrary $\theta>0$,
$P_{\alpha,\theta\nu}$ is a stationary distribution of the
$\mathcal{A}_{\alpha,\theta\nu}$-process. For the same reason as in the
proof of Proposition~\ref{pr3.4} [cf.~(\ref{3.12})], it is sufficient
to prove that
%
\begin{equation}\label{3.21}
\int_{\mathcal{M}_1(E)} {P}_{\alpha,\theta\nu}(d\mu)\mathcal{A}_{\alpha,\theta\nu}
\Phi(\mu)=0
\end{equation}
for $\Phi$ of the form $\Phi(\mu)=\langle\mu,f\rangle^n$ with $f\in
C(E)$ and $n$ being a positive integer. Since Proposition~\ref{pr3.4}
and Lemma~\ref{le3.5} together imply that (\ref{3.21}) holds true for
any $\theta>1$, it is enough to show that the left-hand side of
(\ref{3.21}) defines a real analytic function of $\theta>0$. We claim
that
%
\begin{eqnarray}\label{3.22}
\hspace*{27pt}&& \mathcal{A}_{\alpha,\theta\nu}\Phi(\mu)\nonumber
\\
&&\qquad =  \frac{1}{\Gamma(n)} \sum
_{k=2}^n \pmatrix{n\cr k} (1-\alpha)_{k-2} (\alpha+1)_{n-k} \bigl( \bigl\langle
\mu,f^k \bigr\rangle\langle\mu,f\rangle^{n-k}-\langle\mu,f\rangle^{n} \bigr)\nonumber
\\
&&\quad\qquad{} + \frac{\theta}{(\alpha+1)\Gamma(n)}\nonumber
\\
&&\hphantom{\quad\qquad{} +}{}\times \sum_{k=1}^n \pmatrix{n
\cr k} (1-\alpha)_{k-1} (\alpha)_{n-k} \bigl( \bigl\langle
\nu,f^k \bigr\rangle\langle\mu,f\rangle^{n-k}-\langle\mu,f\rangle
^{n} \bigr)
\\
&&\qquad =  \frac{1}{\Gamma(n)} \sum_{k=2}^n \pmatrix{n
\cr k} (1-\alpha)_{k-2} (\alpha+1)_{n-k} \bigl\langle
\mu,f^k \bigr\rangle\langle\mu,f\rangle^{n-k}\nonumber
\\
&&\quad\qquad{} + \frac{\theta}{(\alpha+1)\Gamma(n)} \sum_{k=1}^n \pmatrix{n
\cr k} (1-\alpha)_{k-1} (\alpha)_{n-k} \bigl\langle
\nu,f^k \bigr\rangle\langle\mu,f\rangle^{n-k}\nonumber
\\
&&\quad\qquad{} -\frac{(\alpha+1)_{n-1}}{(\alpha+1)\Gamma(n)} (\theta+n-1)\langle
\mu,f\rangle^{n}.\nonumber
\end{eqnarray}
The first equality is a special case of (\ref{3.1a}),
and the second one can be shown with the help of Leibniz's formula
\[
(\phi_1\phi_2)^{(n)}(0)=\sum
_{k=0}^n \pmatrix{n\cr k} \phi_1^{(n-k)}(0)
\phi_2^{(k)}(0)
\]
for $\phi_1(t)=(1-t)^{-a}$ and $\phi_2(t)=(1-t)^{-b}$ with
$(a,b)=(\alpha+1,-\alpha-1)$ or $(a,b)=(\alpha,-\alpha)$. In view of
(\ref{3.22}), it is clear that the proof reduces to verifying real
analyticity of $\int
{P}_{\alpha,\theta\nu}(d\mu)\langle\mu,f_1\rangle\cdots
\langle\mu,f_n\rangle$ in $\theta$ for arbitrary $f_1,\ldots,f_n\in
C(E)$.

To this end, we shall exploit the following identity
which is equivalent to~(\ref{3.20}):
%
\begin{equation}\label{3.23}
\qquad\quad\int_{\mathcal{M}_1(E)} {P}_{\alpha,\theta\nu}(d\mu)\langle\mu,1+f\rangle
^{-\alpha} = \int_{\mathcal{M}_1(E)}\mathcal{D}_{\theta\nu}(d\mu)
\bigl\langle\mu, (1+f)^{\alpha} \bigr\rangle^{-1},
\end{equation}
where $f\in B_+(E)$ is arbitrary. Clearly, this remains true for all
bounded Borel functions $f$ on $E$ such that $\inf_{r\in E}f(r)>-1$.
Therefore, for any $t_1,\ldots,t_n\in\mathbf{R}$ with
$|t_1|+\cdots+|t_n|$ being sufficiently small, (\ref{3.23}) for
$f=-\sum_{i=1}^nt_if_i$ is valid, that is,
$I(t_1,\ldots,t_n)=J(t_1,\ldots,t_n)$, where
%
\begin{equation}\label{3.24}
I(t_1,\ldots,t_n) =\int_{\mathcal{M}_1(E)}
{P}_{\alpha,\theta\nu}(d\mu) \Biggl(1- \Biggl\langle\mu,\sum
_{i=1}^nt_if_i \Biggr\rangle\Biggr)^{-\alpha}
\end{equation}
and
%
\begin{equation}\label{3.25}
J(t_1,\ldots,t_n) = \int_{\mathcal{M}_1(E)}
\mathcal{D}_{\theta\nu}(d\mu) \Biggl\langle\mu, \Biggl(1-\sum
_{i=1}^nt_if_i
\Biggr)^{\alpha} \Biggr\rangle^{-1}.
\end{equation}
Noting that $(1-t)^{-\alpha}=1+\sum_{k=1}^{\infty}(\alpha)_kt^k/k!$ as
long as $|t|$ is small enough, we see from (\ref{3.24}) that the
coefficient of the monomial $t_1\cdots t_n$ in the expansion of
$I(t_1,\ldots,t_n)$ is given by
%
\begin{equation}\label{3.26}
(\alpha)_n\int_{\mathcal{M}_1(E)} {P}_{\alpha,\theta\nu}(d\mu)
\langle\mu,f_1\rangle\cdots\langle\mu,f_n\rangle.
\end{equation}
To find the corresponding coefficient for $J(t_1,\ldots,t_n)$,
define
\[
h_{\alpha}(t)=1-(1-t)^{\alpha}= \alpha\sum
_{l=1}^{\infty}(1-\alpha)_{l-1}t^l/l!
\]
and observe from (\ref{3.25}) that
$J(t_1,\ldots,t_n)$ equals
\begin{eqnarray*}
&& \int_{\mathcal{M}_1(E)}\mathcal{D}_{\theta\nu}(d\mu) \Biggl
\langle\mu, 1-h_{\alpha} \Biggl(\sum_{i=1}^nt_if_i
\Biggr) \Biggr\rangle^{-1}
\\
&&\quad =  1+ \sum_{k=1}^{\infty} \int
_{\mathcal{M}_1(E)}\mathcal{D}_{\theta\nu}(d\mu) \Biggl\langle\mu,
h_{\alpha} \Biggl(\sum_{i=1}^nt_if_i
\Biggr) \Biggr\rangle^{k}
\\
&&\quad =  1+ \sum_{k=1}^{\infty}
\alpha^k \int_{\mathcal{M}_1(E)}\mathcal{D}_{\theta\nu}(d
\mu) \sum_{l_1,\ldots,l_k=1}^{\infty} \prod
_{j=1}^k \Biggl\{\frac{(1-\alpha)_{l_j-1}}{l_j!} \Biggl\langle
\mu, \Biggl(\sum_{i=1}^nt_if_i
\Biggr)^{l_j} \Biggr\rangle\Biggr\}.
\end{eqnarray*}
%
One can see that
the coefficient of the monomial $t_1\cdots t_n$
in the expansion of $J(t_1,\ldots,t_n)$ can be expressed as
%
\begin{equation}\label{3.27}
\qquad\sum_{k=1}^n\alpha^k k!\sum
_{\gamma\in\pi(n,k)} \int_{\mathcal{M}_1(E)}
\mathcal{D}_{\theta\nu}(d\mu) \prod_{j=1}^k
\biggl\{\frac{(1-\alpha)_{|\gamma_j|-1}}{|\gamma_j|!} \biggl\langle\mu,\prod_{i\in\gamma_j}f_i
\biggr\rangle\biggr\},
\end{equation}
where $\pi(n,k)$ is the set of partitions $\gamma$ of $\{1,\ldots,n\}$
into $k$ unordered nonempty subsets $\gamma_1,\ldots,\gamma_k$. By
Lemma~2.2 of \cite{E} (or equivalently by Lemma 2.4 of \cite{EG}), each
integral in the above sum is a real analytic function of $\theta>0$.
Hence, so is the integral in (\ref{3.26}) and the stationarity of
$P_{\alpha,\theta\nu}$ with respect to $\mathcal{A}_{\alpha,\theta\nu}$
follows.

It remains to prove the uniqueness of stationary distribution $P$ of
the \mbox{$\mathcal{A}_{\alpha,\theta\nu}$-}pro\-cess for each $\theta>0$. But
this is an immediate consequence of (\ref{3.21}) with $P$ in place of
$P_{\alpha,\theta\nu}$ and (\ref{3.22}), which together determine
uniquely $\int P(d\mu)\langle\mu,f\rangle^n$ and hence the $n$th moment
measure
\[
M_n(dr_1\cdots dr_n): =\int
_{\mathcal{M}_1(E)} P(d\mu)\mu(dr_1)\cdots\mu(dr_n)
\]
for any $n=1,2,\ldots.$ This completes the proof of
Theorem~\ref{th3.2}.
\end{pf*}

It is not clear whether we can derive from (\ref{3.27}) an extension of
the Ewens sampling formula in some explicit and informative form. (See
Remarks after the proof of Lemma~2.2 in \cite{E}.) In view of
(\ref{3.18}), one might think that Pitman's sampling formula would be
applicable. But it is not the case since $\mathcal{D}_m(\mu\mbox{ is
discrete})=1$. The expression (\ref{3.11}) might be rather useful for
such a purpose.

\section{Irreversibility}\label{sec4}
In this section, we discuss reversibility of our processes. In contrast
with the Fleming--Viot diffusion case, we guess that for any
$0<\alpha<1$ and nondegenerate $m$ the $\mathcal{A}_{\alpha,m}$-process
would be irreversible.
Unfortunately, the following result does not give
an affirmative answer in all cases.
However, this does not suggest any possibility of the reversibility
in the exceptional case, which is believed to be
dealt with a different choice of test functions.

%
\begin{theorem}\label{th4.1}
Let $m\in\mathcal{M}(E)^{\circ}$ be given.
Assume that either of the following two conditions holds.
\begin{longlist}[(ii)]
\item[(i)] The support of $m$ has at least three distinct points.

\item[(ii)] The support of $m$ has exactly two points, say $r_1$ and
    $r_2$ and $m(\{r_1\})\ne m(\{r_2\})$.
\end{longlist}
Then the stationary distribution ${P}_{\alpha,m}$ of the
$\mathcal{A}_{\alpha,m}$-process is not a reversible distribution of
it.
\end{theorem}

\begin{pf}
As in the proof of Theorem~\ref{th3.2}, we write $\theta\nu$ instead of
$m$. Thus, $\theta>0$ and $\nu\in\mathcal{M}_1(E)$. Recall that an
equivalent condition to the reversibility of $P_{\alpha,\theta\nu}$
with respect to $\mathcal{A}_{\alpha,\theta\nu}$ is the symmetry
\[
E \bigl[\Phi\mathcal{A}_{\alpha,\theta\nu}\Phi' \bigr] = E \bigl[
\Phi'\mathcal{A}_{\alpha,\theta\nu}\Phi\bigr],\qquad\Phi,
\Phi'\in\mathcal{F}_0,
\]
in which $E[\cdot]$ stands for the expectation with respect to
$P_{\alpha,\theta\nu}$. (See the proof of Theorem 2.3 in \cite{E}.) In
the rest of the proof, we suppress the suffix ``$\alpha,\theta\nu$''
for simplicity. Let $f\in C(E)$ be given and define
$\Phi_n(\mu)=\langle\mu,f\rangle^n$ for each positive integer $n$. We
are going to calculate
%
\begin{equation}\label{4.1}
\Delta:= E [\Phi_2\mathcal{A}\Phi_1 ]-E [
\Phi_1\mathcal{A}\Phi_2 ].
\end{equation}
For this purpose, observe from (\ref{3.22}) that
%
\begin{equation} \label{4.2}
\mathcal{A}\Phi_1(\mu) =\frac{\theta}{\alpha+1} \bigl(\langle\nu,f
\rangle-\langle\mu,f\rangle\bigr),
\end{equation}
%
\begin{eqnarray}\label{4.3}
\mathcal{A}\Phi_2(\mu) 
&=&  \bigl\langle\mu,f^2 \bigr\rangle+\frac{2\alpha\theta}{\alpha+1}\langle\nu,f\rangle\langle\mu,f\rangle
\nonumber\\[-8pt]\\[-8pt]
&&{} + \frac{(1-\alpha)\theta}{\alpha+1} \bigl\langle\nu,f^2 \bigr\rangle-(\theta+1)
\langle\mu,f\rangle^2\nonumber
\end{eqnarray}
and
%
\begin{eqnarray}\label{4.4}
&& \Gamma(3)\mathcal{A}\Phi_3(\mu)\nonumber
\\
&&\qquad =  3(\alpha+1) \bigl\langle
\mu,f^2 \bigr\rangle\langle\mu,f\rangle+(1-\alpha) \bigl\langle
\mu,f^3 \bigr\rangle\nonumber
\\
&&\quad\qquad{} +\frac{\theta}{\alpha+1}\cdot3\alpha(\alpha+1) \langle\nu,f\rangle
\langle\mu,f\rangle^2
\nonumber\\[-12pt]\\[-6pt]
&&\quad\qquad{} +\frac{\theta}{\alpha+1}\cdot3(1-\alpha)\alpha\bigl
\langle\nu,f^2 \bigr\rangle\langle\mu,f\rangle\nonumber
\\
&&\quad\qquad{} +\frac{\theta}{\alpha+1}\cdot(1-\alpha) (2-\alpha) \bigl\langle\nu,f^3 \bigr\rangle-(\alpha+2) (\theta+2)\langle\mu,f\rangle
^3.\nonumber
\end{eqnarray}
Combining (\ref{4.2}) with the stationarity $E[\mathcal{A}\Phi _1]=0$,
we get $E[\langle\mu,f\rangle]=\langle\nu,f\rangle$. Therefore, it is
possible to deduce from (\ref{4.3}) and $E[\mathcal{A}\Phi_2]=0$
\[
(\theta+1)E \bigl[\langle\mu,f\rangle^2 \bigr] =
\frac{2\alpha\theta}{\alpha+1} \langle\nu,f\rangle^2 + \biggl(1+
\frac{(1-\alpha)}{\alpha+1}\theta\biggr) \bigl\langle\nu,f^2 \bigr
\rangle.
\]
Moreover, this equality between quadratic forms is enough to imply
the one between symmetric bilinear forms:
%
\begin{eqnarray}\label{4.5}
&& (\theta+1)E \bigl[\langle\mu,f\rangle\langle\mu,g\rangle\bigr]
\nonumber\\[-8pt]\\[-8pt]
&&\qquad = \frac{2\alpha\theta}{\alpha+1}\langle\nu,f\rangle\langle\nu,g\rangle
+ \biggl(1+\frac{(1-\alpha)}{\alpha+1}\theta\biggr)\langle\nu,fg\rangle,\nonumber
\end{eqnarray}
where $g\in C(E)$ is also arbitrary.
In the rest of the proof, we assume that \mbox{$\langle\nu,f\rangle=0$}.
This makes the calculations below considerably simple.
By (\ref{4.5})
%
\begin{equation}\label{4.6}
M_{1,2}:=E \bigl[\langle\mu,f\rangle\bigl\langle\mu,f^2
\bigr\rangle\bigr] = \frac{(\alpha+1)+(1-\alpha)\theta}{(\alpha
+1)(\theta+1)} \bigl\langle\nu,f^3 \bigr
\rangle.
\end{equation}
The equality $E[\mathcal{A}\Phi_3]=0$ together with (\ref{4.4})
implies that
%
\begin{eqnarray}\label{4.7}
&& (\alpha+2) (\theta+2)E \bigl[\langle\mu,f\rangle^3 \bigr]
\nonumber\\[-9pt]\\[-9pt]
&&\qquad =3(\alpha+1)M_{1,2} +(1-\alpha) \biggl(1+\frac{2-\alpha}{\alpha+1} \theta
\biggr) \bigl\langle\nu,f^3 \bigr\rangle.\nonumber
\end{eqnarray}
These preliminaries help us calculate $\Delta$ in (\ref{4.1}) as follows.
By (\ref{4.3}) and (\ref{4.4})
\begin{eqnarray*}
\Delta &=&  E \biggl[\langle\mu,f\rangle^2 \biggl(-
\frac{\theta}{\alpha+1}\langle\mu,f\rangle\biggr) \biggr] -E \bigl
[\langle\mu,f\rangle\bigl( \bigl\langle\mu,f^2 \bigr\rangle-(\theta+1) \langle
\mu,f\rangle^2\bigr) \bigr]
\\
&=&  \frac{(\alpha+1)+\alpha\theta}{\alpha+1}E \bigl[\langle\mu,f\rangle^3 \bigr]-M_{1,2}
\end{eqnarray*}
and hence (\ref{4.7}) yields
\begin{eqnarray*}
&& (\alpha+1) (\alpha+2) (\theta+2)\Delta
\\
&& \qquad =  \bigl[(\alpha+1)+\alpha\theta\bigr] \biggl[3(\alpha+1)M_{1,2}
+(1-\alpha) \biggl(1+\frac{2-\alpha}{\alpha+1}\theta\biggr) \bigl
\langle\nu,f^3 \bigr\rangle\biggr]
\\
&&\quad\qquad{} -(\alpha+1) (\alpha+2) (\theta+2)M_{1,2}
\\
&&\qquad =  (\alpha+1) (\alpha-1) (2\theta+1)M_{1,2}
\\
&&\quad\qquad{}+ \bigl[(\alpha+1)+
\alpha\theta\bigr](1-\alpha) \biggl(1+\frac{2-\alpha}{\alpha+1}\theta
\biggr) \bigl\langle\nu,f^3 \bigr\rangle.
\end{eqnarray*}
Plugging (\ref{4.6}) into this expression, we obtain
\[
(\alpha+1) (\alpha+2) (\theta+2)\Delta= \frac{1-\alpha}{(\alpha
+1)(\theta+1)} U(\alpha,\theta)
\bigl\langle\nu,f^3 \bigr\rangle,
\]
where
\begin{eqnarray*}
U(\alpha,\theta) &=&  -(\alpha+1) (2\theta+1) \bigl[(\alpha
+1)+(1-\alpha)\theta\bigr]
\\
&&{} + \bigl[(\alpha+1)+\alpha\theta\bigr](\theta+1) \bigl[(\alpha
+1)+(2-\alpha)
\theta\bigr]
\\
&=&  \alpha\theta^2 \bigl[(\alpha+4)+(2-\alpha)\theta\bigr] =: V(
\alpha,\theta).
\end{eqnarray*}
[The second equality between quadratic functions of $\alpha$ is
verified by checking that
$U(-1,\theta)=-3\theta^2(\theta+1)=V(-1,\theta)$,
$U(0,\theta)=0=V(0,\theta)$ and
$U(1,\theta)=\theta^2(\theta+5)=V(1,\theta)$.] Consequently, whenever
$\langle\nu,f\rangle=0$, we have
\[
\Delta=\frac{\alpha(1-\alpha)
\theta^2 [(\alpha+4)+(2-\alpha)\theta]} {
(\alpha+1)^2(\alpha+2)(\theta+1)(\theta+2)} \bigl\langle\nu,f^3 \bigr
\rangle.
\]

Thus, all that remains is to construct an $f\in C(E)$ such that
$\langle\nu,f\rangle=0$ and $\langle\nu,f^3\rangle>0$. Because of the
assumption, we can choose a closed subset $E_0$ of $E$ such that
$0<\nu(E_0)<1/2$. Indeed, in the case (ii) this is trivial while in the
case (i) there exist disjoint closed subsets $E_1,E_2$ and $E_3$ of $E$
such that $\nu(E_1)\nu(E_2)\nu(E_3)>0$ and so $0<\nu(E_i)<1/2$ for some
$i\in\{1,2,3\}$. Letting $g$ denote the indicator function of $E_0$, we
observe that
\begin{eqnarray*}
\bigl\langle\nu, \bigl(g-\langle\nu,g\rangle\bigr)^3 \bigr
\rangle &=&  \bigl\langle\nu,g^3 \bigr\rangle-3 \bigl\langle
\nu,g^2 \bigr\rangle\langle\nu,g\rangle+3\langle\nu,g\rangle
\langle\nu,g\rangle^2-\langle\nu,g\rangle^3
\\
&=&  \nu(E_0)-3\nu(E_0)^2+2
\nu(E_0)^3
\\
&=&  \nu(E_0) \bigl(1-\nu(E_0) \bigr) \bigl(1-2
\nu(E_0) \bigr)>0.
\end{eqnarray*}
Finally, the required $f$ exists since $g$ can be approximated
boundedly and pointwise by a sequence of functions in $C(E)$. The proof
of the theorem is complete.
\end{pf}

It is worth noting that the exceptional case of Theorem~\ref{th4.1}
corresponds to a subclass of the one-dimensional case discussed in
Section~\ref{sec1}, more specifically, the process generated by
(\ref{1.3}) with $c_1=c_2$. There is no reason why this class should be
so special with respect to the reversibility, and it seems that such a
``spatial symmetry'' makes it more subtle to see the asymmetry in time.
The actual difficulty in showing the irreversibility for these
processes along similar lines to the above proof is
that expressions of $E[\Phi_{n_1}\mathcal{A}\Phi_{n_2}]$
with $n_1+n_2\ge4$ as functions of $\alpha$~and~$\theta$
are too complicated to handle.

\section*{Acknowledgment}
The author would like to thank the referees for pointing out some
mistakes and helpful comments on the earlier version of the manuscript.
He is indebted to Professor S. Mano for the reference to \cite{BB}.
Part of the work was carried out under
the ISM Cooperative Research Program
(2012${}\cdot{}$ISM${}\cdot{}$CRP-5008).



\printaddresses

\end{document}